\documentclass[10pt]{amsart}
      \usepackage[mathscr]{eucal}
      \usepackage{amsmath,amsfonts}
      \parskip\smallskipamount
          \newtheorem{theorem}{Theorem}[section]
      
      \newtheorem{proposition}[theorem]{Proposition}
      \newtheorem{corollary}[theorem]{Corollary}
      \newtheorem{lemma}[theorem]{Lemma}

      \makeatletter
      \@addtoreset{equation}{section}
      \makeatother

      \newcommand{\BB}{{\mathbb B}}
      \newcommand{\CC}{{\mathbb C}}
      \newcommand{\NN}{{\mathbb N}}

      \newcommand{\DD}{{\mathbb D}}
      
      \newcommand{\FF}{{\mathbb F}}
      \newcommand{\TT}{{\mathbb T}}

      \newcommand{\cA}{{\mathcal A}}
      
      \newcommand{\cC}{{\mathcal C}}
      
      \newcommand{\cE}{{\mathcal E}}

      \newcommand{\cH}{{\mathcal H}}
      \newcommand{\cK}{{\mathcal K}}
      \newcommand{\cL}{{\mathcal L}}
      
      \newcommand{\cN}{{\mathcal N}}
      \newcommand{\cP}{{\mathcal P}}
      \newcommand{\cR}{{\mathcal R}}

      \newdimen\expt
      \def\boxit#1{\setbox0\hbox{$\displaystyle{#1}$}
            \hbox{\lower.4\expt
       \hbox{\lower3\expt\hbox{\lower\dp0
            \hbox{\vbox{\hrule height.4\expt
       \hbox{\vrule width.4\expt\hskip3\expt
            \vbox{\vskip3\expt\box0\vskip2\expt}%
       \hskip3\expt\vrule width.4\expt}\hrule height.4\expt}}}}}}
      
\begin{document}
       \pagestyle{myheadings}
      \markboth{ Gelu Popescu}{   Free  holomorphic functions
        and  interpolation }

      \title [   Free  holomorphic functions
        and interpolation ]
      {    Free  holomorphic functions
        and  interpolation
      }
        \author{Gelu Popescu}
\date{September 14, 2007}
      \thanks{Research supported in part by an NSF grant}

      \address{Department of Mathematics, The University of Texas
      at San Antonio \\ San Antonio, TX 78249, USA}
      \email{\tt gelu.popescu@utsa.edu}

\begin{abstract}
In this  paper we  obtain a noncommutative multivariable analogue of
the classical Nevanlinna-Pick interpolation problem for analytic
functions with positive real parts on the open unit disc. Given a
function $f:\Lambda\to \CC$, where $\Lambda$ is  an arbitrary subset
of the open unit ball $\BB_n:=\{z\in \CC^n: \|z\|<1\}$, we find
necessary and sufficient conditions for the existence  of  a free
holomorphic function $g$ with complex
 coefficients on the noncommutative open unit ball $[B(\cH)^n]_1$
such that
$$
 \ \text{\rm Re} \  g \geq 0 \ \text{ and } \ g(z)=f(z),\
z\in \Lambda, $$ where $B(\cH)$ is the algebra of all bounded linear
operators on a Hilbert space $\cH$.
 The proof employs  several
results from noncommutative multivariable operator theory and
 a noncommutative Cayley transform  (introduced and studied in the present paper) acting from
the set of all free holomorphic functions with positive real parts
to the set of all bounded free holomorphic functions. All  the
results of this paper are obtained in the more general setting of
free holomorphic functions with operator-valued coefficients. As
consequences, we deduce some  results concerning operator-valued
analytic
 interpolation on the unit ball $\BB_n$.
\end{abstract}

      \maketitle

\bigskip

\bigskip

\section*{Introduction}

The classical Nevanlinna-Pick interpolation problem (\cite{Pic},
\cite{N})  for analytic functions with positive real parts on the
open unit disc $\DD:=\{\lambda\in \CC:\ |\lambda|<1\}$ is the
following: given $m$ distinct points $\lambda_1,\ldots, \lambda_m$
in $\DD$ and $m$ complex numbers $w_1,\ldots, w_m$, find an analytic
function in the open unit disc  with $\text{\rm Re}\, f(\lambda)\geq
0$, $\lambda\in \DD$, such that $f(\lambda_i)=w_i$ for any
$i=1,\ldots,m$. It was proved that this interpolation problem has a
solution if and only if the matrix
$$
\left[ \begin{matrix}
\frac{w_1+\overline{w}_1}{1-\lambda_1\overline{ \lambda}_1}&\cdots&
\frac{w_1+\overline{w}_m}{1-\lambda_1\overline{ \lambda}_m}\\
\vdots &\vdots& \vdots\\
\frac{w_m+\overline{w}_1}{1-\lambda_m\overline{ \lambda}_1}&\cdots &
\frac{w_m+\overline{w}_m}{1-\lambda_m\overline{ \lambda}_m}
\end{matrix}
\right]
$$
is positive semidefinite. In this case, $f$ has the integral
Riesz-Herglotz  (\cite{Ri}, \cite{Her}) representation
$$
f(\lambda)=\int_{\TT}\frac{e^{i\theta}+\lambda}{e^{i\theta}-\lambda}
d \mu(e^{i\theta}) +i\text{\rm Im}\, f(0), \quad \lambda\in \DD,
$$
for some finite positive Borel measure $\mu$ on the unit circle
$\TT:=\{\lambda\in \CC:\ |\lambda|=1\}$. The Nevanlinna-Pick
interpolation problem has been studied further by several authors
and generalized to various settings (see \cite{SzK}, \cite{S},
\cite{AMc2},  \cite{FF-book}, \cite{Po-interpo}, \cite{ArPo2},
\cite{DP},   \cite{EP},   \cite{Po-entropy}, \cite{MuSo2},  and the
references there in).

In the last two decades, significant progress  has been made in
noncommutative multivariable operator theory regarding
 noncommutative dilation theory and its applications to
    interpolation in several variables (\cite{Po-isometric},
     \cite{Po-intert},
   \cite{Po-unitary},  \cite{Po-varieties},
    \cite{Po-analytic},
\cite{BTV}, \cite{Po-nehari}), and unitary invariants for $n$-tuples
of operators (\cite{Po-charact},   \cite{Arv2}, \cite{Po-curvature},
   \cite{Po-similarity}, \cite{BT},
   \cite{Po-unitary}). In related areas of research,
  we remark
 the work  of Helton,  McCullough, Putinar, and  Vinnikov, on symmetric noncommutative
  polynomials (\cite{He}, \cite{He-C},
 \cite{He-C-P2}, \cite{He-C-V}), and the work of Muhly and Solel  on
 representations of tensor algebras over
 $C^*$ correspondences (see \cite{MuSo1}, \cite{MuSo2}).

In \cite{Po-holomorphic}, \cite{Po-pluriharmonic} we developed a
theory of holomorphic (resp. pluriharmonic)  functions in several
noncommuting (free) variables and provide a framework for  the study
of arbitrary
 $n$-tuples of operators on a Hilbert space $\cH$. Several classical
 results from complex analysis have
 free analogues in
 this  noncommutative multivariable setting.
This theory enhances our program to develop a {\it free}
  analogue of
  Sz.-Nagy--Foia\c s theory \cite{SzF-book}, for row contractions.
We recall that a map $f:[B(\cH)^n]_1\to B(\cH)$ is called free
holomorphic function  with scalar coefficients on the noncommutative
open  unit ball
$$
[B(\cH)^n]_1:=\{(X_1,\ldots, X_n)\in B(\cH)^n: \ \|X_1X_1^*+\cdots +
X_nX_n^*\|^{1/2}<1\}
$$
if $f$ has a representation
$$
f(X_1,\ldots, X_n)=\sum_{k=0}^\infty \sum_{|\alpha|=k}
 a_\alpha X_\alpha,\qquad (X_1,\ldots, X_n)\in [B(\cH)^n]_1,
$$
where $\{a_\alpha\}_{\alpha\in \FF_n^+}$ are complex numbers  with \
$\limsup\limits_{k\to\infty} (\sum\limits_{|\alpha|=k}
|a_\alpha|^2)^{1/2k}\leq 1$, and $\FF_n^+$ is the free semigroup
with $n$ generators. In the particular case when $n=1$ and
$\cH=\CC$, we recover the analytic functions on the open unit disc.

The main goal of the present paper is to obtain an analogue of the
classical Nevanlinna-Pick interpolation problem, for free
holomorphic functions $f$ with  $\text{\rm Re}\, f\geq 0$, i.e.,
$$\text{\rm Re
}f(X_1,\ldots, X_n)\,\geq 0,\qquad (X_1,\ldots, X_n)\in
[B(\cH)^n]_1,$$
 for  any Hilbert space
$\cH$. We prove that, given   a function $f:\Lambda\to \CC$, where
$\Lambda$ is  an arbitrary  subset of the open unit ball $\BB_n$,
there exists a free holomorphic function $g$ with complex
 coefficients
such that
$$
 \ \text{\rm Re}\ g \geq 0 \ \text{ and } \ g(z)=f(z),\
z\in \Lambda, $$
 if and only if
    the map
   \begin{equation*}
 \Lambda\times \Lambda \ni
(z,w)\mapsto \frac{f(z)+ \overline{f(w)}}{1-\left<z,w\right>}\in \CC
\end{equation*}
 is positive  semidefinite. In this case, we have a Riesz-Herglotz type representation for
 $g$, i.e.,
 $$
g(X)=(\mu \otimes\text{\rm id})[(I+R_X)(I-R_X)^{-1}]+ i\text{\rm
Im}\, g(0)
$$
 for any $ X:=(X_1,\ldots, X_n)\in [B(\cH)^n]_1$,
 where
$R_X:=R_1^*\otimes X_1  +\cdots + R_n^*\otimes X_n $ is the
reconstruction operator  and $\mu$ is a completely positive linear
functional on the Cuntz-Toeplitz algebra $C^*(R_1,\ldots,R_n)$
generated by the right creation operators on the full Fock space
with $n$ generators. We should add that all the results of this
paper are obtained in the more general setting of free holomorphic
functions with operator-valued coefficients.

 As consequences, we deduce some  results concerning  operator-valued analytic
 interpolation on the unit ball $\BB_n$. In particular, in the scalar case  when  $f:\Lambda\to
 \CC$ satisfies the above-mentioned positivity
 condition,  we show that
  $f$  has an analytic extension
   $\varphi:\BB_n\to \CC$   such that
    the map
  $$ \BB_n\times \BB_n\ni
(z,w)\mapsto \frac{\varphi(z)+
\overline{\varphi(w)}}{1-\left<z,w\right>}\in \CC
$$
 is positive  semidefinite, and $\varphi$ has  the  Riesz-Herglotz type representation
$$
\varphi(z)=\nu[(I+z_1R_1^*+\cdots + z_nR_n^*)(I-z_1R_1^*-\cdots -
z_nR_n^*)^{-1}]+ i\text{\rm Im}\, \varphi(0)
$$
for  any $z:=(z_1,\ldots, z_n)\in \BB_n$, where $\nu$ is a
completely positive linear functional on the Cuntz-Toeplitz algebra
$C^*(R_1,\ldots,R_n)$.

In Section 1,  we introduce a noncommutative Cayley transform acting
from the set of all free holomorphic functions with positive real
parts to the set of all bounded free holomorphic functions. This
transform plays a crucial role in this paper and, due to its
properties,  enable  us to
 use several results from noncommutative multivariable operator
theory in order to prove our interpolation results  in Section 2 and
Section 3.

Finally, we mention that the noncommutative Cayley transform   will
be a key tool in a forthcoming paper, where we study free
pluriharmonic majorants and obtain a description of all solutions of
a generalization of the noncommutative commutant lifting theorem
(see \cite{Po-isometric}, \cite{Po-intert}).

\bigskip

\section{Noncommutative Cayley transforms}
We introduce a Cayley type transform acting from  the set of all
 free holomorphic functions
  with positive real parts to the set of all bounded free holomorphic functions.
  This transform will play an important role in the next
  section where we  solve a  Nevanlinna-Pick   type  interpolation problem
  for free holomorphic
  functions with positive real parts.

First, we need   a     result concerning Cayley transforms of
bounded accretive operators  which traces back to  von Neumann
\cite{von1} (see also \cite{SzF-book}). Since we could not find a
precise reference for the following version, we sketch a proof for
completeness.

\begin{proposition} \label{Re}
An operator $T\in B(\cH)$ is a contraction  $( \|T\|\leq 1 )$   such
that $I-T$ is invertible  if and only is  there exists $A\in B(\cH)$
with $\text{\rm Re}\,A\geq 0$ such that $T$ is equal to the Cayley
transform of $A$, i.e.,
$$
T=\cC(A):=(A-I)(A+I)^{-1}.
$$
Moreover,  in this case, $A$ is the inverse Cayley transform of $T$,
i.e., $$A= \cC^{-1}(T):=(I+T)(I-T)^{-1}.
$$
\end{proposition}
\begin{proof}
Assume that  $\|T\|\leq 1$ and $I-T$ is invertible. Define the
operator $A:=(I+T)(I-T)^{-1}$ and notice that
\begin{equation*}
\begin{split}
A +A^* &=(I-T)^{-1}(I+T)+ (I+T^*)(I-T^*)^{-1}\\
&=2(I-T)^{-1}(I-TT^*)(I-T^*)^{-1}\geq 0.
\end{split}
\end{equation*}
Therefore $\text{\rm Re}\, A\geq 0$. On the other hand, we have
$A=2(I-T)^{-1}-I$ which implies that $I+A$ is invertible.
Consequently, one can easily see that $T=(A-I)(A+I)^{-1}$.

Conversely, let $A\in B(\cH)$  be such that  $\text{\rm Re}\,A\geq
0$ and note that
\begin{equation*}
\begin{split}
 \|Ah+h\|^2&=\|Ah\|^2+\left<\text{\rm Re}\, A h, h\right>
+\|h\|^2\\
&\geq \|Ah\|^2-\left<\text{\rm Re}\, A h, h\right>
+\|h\|^2=\|Ah-h\|^2
\end{split}
\end{equation*}
for any $h\in \cH$. Hence, we deduce that $\|Ah+h\|\geq \|h\|$,
$h\in \cH$. Similar inequalities, as above, hold if we replace $A$
with $A^*$. Therefore, we also have  $\|A^*h+h\|\geq \|h\|$, $h\in
\cH$. Since the operators $A+I$ and $A^*+I$ are bounded below, we
infer that $A+I$ is invertible. On the other hand, the
above-mentioned  inequalities, imply that the operator $T:\cH\to
\cH$ defined by $T(A+I)h:=(A-I)h$, $h\in \cH$,  is a contraction and
$T=(A-I)(A+I)^{-1}$. Hence, an easy calculation shows that
$I-T=2(A+I)^{-1}$ and, consequently, $I-T$ is invertible. Moreover,
we  have  $A=(I+T)(I-T)^{-1}$.  The proof is complete.
\end{proof}

We  recall  from \cite{Po-holomorphic} a few facts concerning free
holomorphic functions on the  noncommutative open ball
 $[B(\cH)^n]_1 $.
Let $\FF_n^+$ be the unital free semigroup on $n$ generators
$g_1,\ldots, g_n$ and the identity $g_0$.  The length of $\alpha\in
\FF_n^+$ is defined by $|\alpha|:=0$ if $\alpha=g_0$  and
$|\alpha|:=k$ if
 $\alpha=g_{i_1}\cdots g_{i_k}$, where $i_1,\ldots, i_k\in \{1,\ldots, n\}$.
If $(X_1,\ldots, X_n)\in B(\cH)^n$, where $B(\cH)$ is the algebra of
all bounded linear operators on the Hilbert space $\cH$,    we
denote $X_\alpha:= X_{i_1}\cdots X_{i_k}$  and $X_{g_0}:=I_\cH$.

 Let
$\{A_{(\alpha)}\}_{\alpha\in \FF_n^+}$ be a sequence
 of bounded linear operators  on a Hilbert space   $\cE$ and
 define
 $R\in [0,\infty]$ by setting
$$ \frac {1} {R}:= \limsup_{k\to\infty} \left\|\sum_{|\alpha|=k}
A_{(\alpha)}^* A_{(\alpha)}\right\|^{\frac{1} {2k}}. $$
We call   $R$
 {\it radius of
convergence} of the formal power series $\sum_{\alpha\in \FF_n^+}
A_{(\alpha)} \otimes  Z_\alpha $
 in  noncommuting indeterminates
  $Z_1,\ldots, Z_n$,
where $Z_\alpha:=Z_{i_1}\cdots Z_{i_k}$ if $\alpha=g_{i_1}\cdots
g_{i_k}$
 and $Z_{g_0}:=1$.

A map $F:[B(\cH)^n]_1\to B(\cE)\otimes_{min} B( \cH)$ is called a
  {\it free
holomorphic function} on  $[B(\cH)^n]_1$  with coefficients in
$B(\cE)$ if there exists $A_{(\alpha)}\in B(\cE)$, $\alpha\in
\FF_n^+$, such that the formal power series $\sum_{\alpha\in
\FF_n^+} A_{(\alpha)} \otimes  Z_\alpha$ has radius of convergence
$\geq 1$ and such that
$$
F(X_1,\ldots, X_n)=\sum\limits_{k=0}^\infty
\sum\limits_{|\alpha|=k}A_{(\alpha)}\otimes  X_\alpha,\quad
(X_1,\ldots, X_n)\in [B(\cH)^n]_1,
$$
where the series converges in the operator  norm topology.  We
proved in \cite{Po-holomorphic}
 that the following statements are equivalent:
 \begin{enumerate}
 \item[(i)] the series $\sum\limits_{k=0}^\infty
\sum\limits_{|\alpha|=k}A_{(\alpha)}\otimes X_\alpha$ is convergent
in the operator norm for any $(X_1,\ldots, X_n)\in [B(\cH)^n]_1$ and
any Hilbert space $\cH$;
\item[(ii)] $\limsup\limits_{k\to\infty}\left\|\sum\limits_{|\alpha|=k}
A_{(\alpha)}^* A_{(\alpha)}\right\|^{1/2k}\leq  1 $;
\item[(iii)]  the series
 $\sum_{k=1}^\infty \sum_{|\alpha|=k} A_{(\alpha )}\otimes r^{|\alpha|}S_\alpha$
 is convergent in the operator norm  for any $r\in [0,1)$.
\end{enumerate}
We remark that if $f$ is a free holomorphic function on
$[B(\cH)^n]_1$, then its scalar representation (when $\cH=\CC$) is
an operator-valued analytic function $z\mapsto f(z)$ on the unit
ball $\BB_n$. The converse is not true in general if   $n\geq 2$.
The set of all free holomorphic functions on $[B(\cH)^n]_1$ with
coefficients in $B(\cE)$ is denoted by $Hol(B(\cH)^n_1)$.
Let $H^\infty(B(\cH)^n_1)$  denote the set of  all elements $F$ in
$Hol(B(\cH)^n_1)$  such that
$$
\|F\|_\infty:=\sup  \|F(X_1,\ldots, X_n)\|<\infty,
$$
where the supremum is taken over all $n$-tuples  of operators
$(X_1,\ldots, X_n)\in [B(\cH)^n]_1$ and any Hilbert space $\cH$.
According to \cite{Po-holomorphic} and \cite{Po-pluriharmonic},
$H^\infty(B(\cH)^n_1)$ can be identified to the operator algebra
$B(\cE)\bar\otimes F_n^\infty$ (the weakly closed algebra generated
by the spatial tensor product), where $F_n^\infty$ is the
noncommutative analytic Toeplitz algebra (see \cite{Po-von},
\cite{Po-multi}, \cite{Po-analytic}).

 We say that $G$ is a  self-adjoint  {\it free pluriharmonic
function} on $[B(\cH)^n]_1$ if there exists a free holomorphic
function $F$ on $[B(\cH)^n]_1$ such that $G=\text{\rm Re}\, F$,
i.e.,
$$
G(X_1,\ldots, X_n)=\text{\rm Re}\, F(X_1,\ldots,
X_n):=\frac{1}{2}(F(X_1,\ldots, X_n)+ F(X_1,\ldots, X_n)^*)
$$
for  $(X_1,\ldots, X_n)\in [B(\cH)^n]_1$.   An arbitrary free
pluriharmonic function on $[B(\cH)^n]_1$  has the form
 $H:=H_1+iH_2$, where $H_1$ and $H_2$ are self-adjoint
free harmonic functions on $[B(\cH)^n]_1$.  We remark that in the
particular case when $n=1$, a function is free pluriharmonic on
$[B(\cH)]_1$ if and only if its scalar representation (when
$\cH=\CC$) is harmonic on the open unit disc $\DD$. The study of
free pluriharmonic functions was pursued in our extensive paper
\cite{Po-holomorphic}. Some of those results will be used in the
present paper.

In what follows we consider a few preliminary results concerning
formal power series   in noncommutative indeterminates and
operator-valued coefficients.

\begin{lemma}
\label{inverse}
 Let $f=\sum_{\alpha\in \FF_n^+} A_{(\alpha)}\otimes
Z_\alpha$ be a formal power series in noncommutative indeterminates
$Z_1,\ldots, Z_n$, coefficients  $A_{(\alpha)}$  in $B(\cE)$, and
such that $A_{(g_0)}$ is an invertible operator.  Then $f$ has an
inverse as a power series with operator-valued coefficients.
\end{lemma}
\begin{proof} Denote $A_{(0)}:=A_{(g_0)}$
 and note
  that $f=(A_{(0)}\otimes 1)(1-g)$, where $g=\sum_{|\alpha|\geq
1} D_{(\alpha)}\otimes Z_\alpha$ and $D_{(\alpha)}:=-A_{(0)}^{-1}
A_{(\alpha)}$. For each $m\in \NN$, $g^m$ defines a power series
$\sum_{|\alpha|\geq m} C_{(\alpha)}\otimes Z_\alpha$ and it makes
sense to consider the formal power series
$$
\varphi=1+g+g^2+\cdots.
$$
We call  the operator $C_{(\alpha)}$ the $\alpha$-coefficient of
$g^m$. Notice that if $m>k$, then the term $g^m$ has all
coefficients of order $\leq k$ equal to $0$. Thus, if $\alpha\in
\FF_n^+$ with $|\alpha|\leq k$, then we may define the
$\alpha$-coefficient of $\varphi$ as the $\alpha$-coefficient of the
finite sum $1+g+ g^2+\cdots + g^m$. Notice  also that $
\varphi=1+\sum_ {|\alpha|\geq 1} F_{(\alpha)}\otimes  Z_\alpha, $
where
\begin{equation*}
 F_{(\alpha)}= \sum_{j=1}^{|\alpha|} \sum_{{\gamma_1\cdots
\gamma_j=\alpha }\atop {|\gamma_1|\geq 1,\ldots, |\gamma_j|\geq 1}}
D_{(\gamma_1)}\cdots D_{(\gamma_j)}   \quad \text{ for } \
|\alpha|\geq 1.
\end{equation*}
Since $(1-g)\varphi=\varphi (1-g)=1$, we have $(1-g)^{-1}=\varphi$.
Now, it is clear that the inverse of $f$ satisfies the relations
$$
f^{-1}=(1-g)^{-1}(A_{(0)}^{-1}\otimes 1)=A_{(0)}^{-1}\otimes 1+
\sum_{|\alpha|\geq 1} F_{(\alpha)}A_{(0)}^{-1}\otimes Z_\alpha.
$$
The proof is complete.
\end{proof}

Let $\widetilde \CC_{+}[Z_1,\ldots, Z_n]$ (resp.~$\widetilde
\CC_{-}[Z_1,\ldots, Z_n]$) be the set of all  formal power series
$f=\sum_{\alpha\in \FF_n^+} A_{(\alpha)}\otimes Z_\alpha$ such that
$I+ A_{(0)}$  (resp.~$I- A_{(0)}$) is an invertible operator in
$B(\cE)$. Denote by $\widetilde \CC_0[Z_1,\ldots, Z_n]$
(resp.~$\widetilde \CC_1[Z_1,\ldots, Z_n]$) the set of all formal
power series in noncommutative indeterminates $Z_1,\ldots, Z_n$,
coefficients in $B(\cE)$,  and constant term $0$ (resp.~$1$).

\begin{proposition}\label{Cayley-series}
The Cayley transform $\widetilde \Gamma: \widetilde
\CC_{+}[Z_1,\ldots, Z_n]\to \widetilde \CC_{-}[Z_1,\ldots, Z_n]$
defined by
$$\widetilde \Gamma(f):=(f-1)(1+f)^{-1},\quad f\in \widetilde
\CC_{+}[Z_1,\ldots, Z_n]
$$
is a bijection, and its inverse satisfies the equation
$$
\widetilde \Gamma^{-1}(g)=(1-g)^{-1}(g+1),\quad g\in \widetilde
\CC_{-}[Z_1,\ldots, Z_n].
$$

In particular, $\widetilde\Gamma$  is a bijection   from $\widetilde
\CC_{1}[Z_1,\ldots, Z_n]$ to $\widetilde \CC_{0}[Z_1,\ldots, Z_n]$.
\end{proposition}

\begin{proof}
Due to  Lemma \ref{inverse}, the map $\widetilde \Gamma$ is
well-defined and the coefficient $B_{(0)}$ of the power series
$\widetilde \Gamma(f)=\sum_{\alpha\in \FF_n^+} B_{(\alpha)}\otimes
Z_\alpha $ satisfies the equation
$B_{(0)}=(A_{(0)}-I)(I+A_{(0)})^{-1}$. Consequently, we have
$I-B_{(0)}=2(I+A_{(0)})^{-1}$, which shows that the power series
$\widetilde \Gamma(f)$ is in $\CC_{-}[Z_1,\ldots, Z_n]$.

If $f_1,f_2\in \widetilde \CC_{+}[Z_1,\ldots, Z_n]$ and  $\widetilde
\Gamma(f_1)=\widetilde \Gamma (f_2)$, then
$(f_1-1)(1+f_2)=(1+f_1)(f_2-1)$, whence $f_1=f_2$. To prove  that
the Cayley transform is surjective, let $g\in \widetilde
\CC_{-}[Z_1,\ldots, Z_n]$.   It is easy to see that, due to Lemma
\ref{inverse},  the power series  $(1-g)^{-1}(g+1)$ is in
$\widetilde \CC_{+}[Z_1,\ldots, Z_n]$ and
\begin{equation*}
\begin{split}
 \widetilde \Gamma [(1-g)^{-1}(g+1)]&=
 \left[(1-g)^{-1}(g+1)-1\right][ 1+(1-g)^{-1}(g+1)]^{-1}\\
 &=\left[(1-g)^{-1}(g+1)-1\right]\left[(1-g)^{-1} (1-g+g+1)\right]^{-1} \\
 &=g.
\end{split}
\end{equation*}
Therefore, $\widetilde \Gamma$ is a bijection. The last part of this
proposition is now obvious.
\end{proof}

Let $H_n$ be an $n$-dimensional complex  Hilbert space with
orthonormal
      basis
      $e_1$, $e_2$, $\dots,e_n$, where $n\in\{1,2,\dots\}$.
       We consider the full Fock space  of $H_n$ defined by
      $$F^2(H_n):=\CC 1\oplus \bigoplus_{k\geq 1} H_n^{\otimes k},$$
      where  $H_n^{\otimes k}$ is the (Hilbert)
      tensor product of $k$ copies of $H_n$.
      Define the left  (resp. right) creation
      operators  $S_i$ (resp.~$R_i$), $i=1,\ldots,n$, acting on $F^2(H_n)$  by
      setting
      $$
       S_i\varphi:=e_i\otimes\varphi, \quad  \varphi\in F^2(H_n),
      $$
       (resp.~$
       R_i\varphi:=\varphi\otimes e_i, \quad  \varphi\in F^2(H_n).
      $)
The noncommutative disc algebra $\cA_n$ (resp.~$\cR_n$) is the norm
closed algebra generated by the left (resp.~right) creation
operators and the identity. The   noncommutative analytic Toeplitz
algebra $F_n^\infty$ (resp.~$R_n^\infty$)
 is the the weakly
closed version of $\cA_n$ (resp.~$\cR_n$). These algebras were
introduced in \cite{Po-von} in connection with a noncommutative von
Neumann inequality (see \cite{von} for the classical case when
$n=1$). They
 have  been studied
    in several papers
\cite{Po-charact},  \cite{Po-multi},  \cite{Po-funct},
\cite{Po-analytic}, \cite{Po-disc}, \cite{Po-poisson},
   and recently in
  \cite{DP1}, \cite{DP2},   \cite{DP},
  \cite{ArPo2}, \cite{Po-curvature},  \cite{DKP},  \cite{PPoS},
   \cite{Po-similarity},     and \cite{Po-unitary}.

Let $\FF_n^+$ be the unital free semigroup on $n$ generators
      $g_1,\dots,g_n$, and the identity $g_0$.
       We denote $e_\alpha:=
e_{i_1}\otimes\cdots \otimes  e_{i_k}$  if $\alpha=g_{i_1}\cdots
g_{i_k}$, where $i_1,\ldots, i_k\in \{1,\ldots,n\}$, and
$e_{g_0}:=1$. Note that $\{e_\alpha\}_{\alpha\in \FF_n^+}$ is an
orthonormal basis for $F^2(H_n)$. Let  $\cP^{(m)}$, $m=0,1,\ldots$,
be the set of all polynomials of degree $\leq m$   in $e_1,\ldots,
e_n$, i.e.,
$$
\cP^{(m)}:=\text{ \rm span} \{ e_\alpha: \ \alpha\in \FF_n^+,
|\alpha|\leq m\},
$$
and define the nilpotent operators $S_i^{(m)}: \cP^{(m)}\to
\cP^{(m)}$ by
$$
S_i^{(m)}:=P_{\cP^{(m)}} S_i |_{\cP^{(m)}},\quad i=1,\ldots, n,
$$
where $S_1,\ldots, S_n$ are the left creation operators on the Fock
space $F^2(H_n)$ and $P_{\cP^{(m)}}$ is the orthogonal projection of
$F^2(H_n)$ onto $\cP^{(m)}$. Notice that $S_\alpha^{(m)}=0$ if
$|\alpha|\geq m+1$. According to \cite{Po-unitary}, the $n$-tuple of
operators $(S_1^{(m)},\ldots, S_n^{(m)})$ is the universal model for
row contractions $(T_1,\ldots, T_n)$  with $T_\alpha=0$ for
$|\alpha|\geq m+1$, and the following constrained von Neumann type
inequality    holds:
\begin{equation*}
 \|p(T_1,\ldots, T_n)\|\leq \|p(S_1^{(m)},\ldots,
S_n^{(m)})\|
\end{equation*}
for any noncommutative polynomial $p(X_1,\ldots,
X_n)=\sum_{|\alpha|\leq k} A_{(\alpha)}\otimes X_\alpha$, \ $k\in
\NN$.
   We
also know that $f\in H^\infty(B(\cH)^n_1)$ if and only if \
 $\sup_{m\in \NN} \|f(S_1^{(m)},\ldots, S_n^{(m)})\|<\infty$.
Moreover, in this case, we have
$$\|f\|_\infty=\sup_{m\in \NN} \|f(S_1^{(m)},\ldots, S_n^{(m)})\|.
$$

We recall from \cite{Po-pluriharmonic}  a few  properties concerning
 positive free pluriharmonic    functions which will be used in the present paper. We say that a free
 pluriharmonic function $g$ is positive if $g(X)\geq 0$ for any
$X:=(X_1,\ldots, X_n)\in [B(\cH)^n_1]$ and any Hilbert space.

 \begin{proposition}\label{properties}
 Let $f$ be  a free holomorphic function on $[B(\cH)^n]_1$  with  the
 representation
 $$
 f(X_1,\ldots, X_n):=\sum_{k=0}^\infty \sum_{|\alpha|=k}
 A_{(\alpha)}\otimes X_\alpha,\quad (X_1,\ldots, X_n)\in
 [B(\cH)^n]_1.
 $$
 Then
 \begin{enumerate}
 \item[(i)]  $\text{\rm Re}\,
f\geq 0$ if and only if \ $\text{\rm Re}\, f(S_1^{(m)},\ldots,
S_n^{(m)})\geq 0$ for any $m\geq 0$;
\item[(ii)]
if \  $\text{\rm Re}\, f\geq 0$, then
$$
\left\|\sum_{|\alpha|=k} A_{(\alpha)}^* A_{(\alpha)}\right\|^{1/2}
\leq \|A_{(0)}+ A_{(0)}^*\|
$$
for any $k\geq 1$.

 \end{enumerate}

 \end{proposition}

Denote by $Hol^+(B(\cH)^n_1)$ the set of all   free holomorphic
functions $f$  with coefficients in $B(\cE)$, where $\cE$ is a
separable Hilbert space,  such that $\text{\rm Re}\, f\geq 0$, and
 let
$H_-^\infty(B(\cH)^n_1)$ be the set of all bounded  free holomorphic
functions $f(X_1,\ldots, X_n)=\sum_{\alpha\in \FF_n^+}
A_{(\alpha)}\otimes X_\alpha$   such that $I_\cE- A_{(0)}$ is an
invertible operator in $B(\cE)$. Consider also the following sets:
\begin{equation*}
\begin{split}
Hol_1^+(B(\cH)^n_1)&:=\left\{ f\in Hol^+(B(\cH)^n_1):\
f(0)=I \right\}, \\
H_0^\infty(B(\cH)^n_1)&:=\left\{ g\in H^\infty(B(\cH)^n_1): \
g(0)=0\right\}.
\end{split}
\end{equation*}

 We introduce the {\it noncommutative Cayley  transform}
$$\Gamma: Hol^+(B(\cH)^n_1)\to
\left[H_-^\infty(B(\cH)^n_1)\right]_{\leq 1}\quad \text{ defined by
} \ \Gamma f:=g,$$
 where $g\in H_-^\infty(B(\cH)^n_1)$ is   uniquely determined by
the formal power series
 $\widetilde\Gamma (\widetilde  f):=( \widetilde
 f-1)(1+\widetilde f)^{-1}$, where
$\widetilde f$ is the power series associated with $f$. Of course,
 we need  to show that  $\Gamma$ is well-defined.

The following result plays a key role in the present paper.

\begin{theorem}
\label{Cayley1} The noncommutative Cayley transform  $\Gamma$ is a
bijection between $Hol^+(B(\cH)^n_1)$ and the unit ball of
$H_-^\infty(B(\cH)^n_1)$. In particular, $\Gamma$ is a bijection
between $Hol_1^+(B(\cH)^n_1)$ and the unit ball of
$H_0^\infty(B(\cH)^n_1) $.

\end{theorem}
\begin{proof}
First, we show that the map $\Gamma$ is well-defined and $\Gamma f$
is in  the unit ball of $H_-^\infty(B(\cH)^n_1)$. Let $f$ be in
$Hol^+(B(\cH)^n_1)$ and have the representation $f(X_1,\ldots,
X_n):=B_{(0)}\otimes I+\sum_{k=1}^\infty \sum_{|\alpha|=k}
B_{(\alpha)}\otimes X_\alpha$. Since $\text{\rm Re}\, f(X)\geq 0$
for any $X \in [B(\cH)^n_1]$, we deduce that $B_{(0)}+B_{(0)}^*\geq
0$. As in the proof of Proposition \ref{Re}, we  can show that
$I_\cE+B_{(0)}$ is invertible. Taking into account that
$S_\alpha^{(m)}=0$ if $|\alpha\geq m+1$, we have
$f(S_1^{(m)},\ldots, S_n^{(m)})=\sum_{|\alpha|\leq m}
B_{(\alpha)}\otimes S_\alpha^{(m)}$.
 Note  that
\begin{equation*}
\begin{split}
H(S_1^{(m)},\ldots, S_n^{(m)})&:=2[f(S_1^{(m)},\ldots, S_n^{(m)})^*
+ f(S_1^{(m)},\ldots, S_n^{(m)})]\\
&=[I+ f(S_1^{(m)},\ldots, S_n^{(m)})^*][I+ f(S_1^{(m)},\ldots,
S_n^{(m)})]\\
&\qquad - [f(S_1^{(m)},\ldots, S_n^{(m)})^*-I] [f(S_1^{(m)},\ldots,
S_n^{(m)})-I].
\end{split}
\end{equation*}
Since $H(S_1^{(m)},\ldots, S_n^{(m)})\geq 0$, we deduce  that
$$
\|[I+ f(S_1^{(m)},\ldots, S_n^{(m)})]x\|\geq \|[f(S_1^{(m)},\ldots,
S_n^{(m)})-I]x\|
$$
for any $x\in\cE\otimes \cP^{(m)}$. Consequently, there exists a
contraction $ G_m:\cE\otimes\cP^{(m)}\to \cE\otimes\cP^{(m)}$ such
that
\begin{equation}
\label{Gm}
 G_m[I+ f(S_1^{(m)},\ldots,
S_n^{(m)})]=f(S_1^{(m)},\ldots, S_n^{(m)})-I.
\end{equation}
Using the fact that    the operator $I_\cE+B_{(0)}$ is invertible,
we have
$$I+f(S_1^{(m)},\ldots, S_n^{(m)})=[(I_\cE+B_{(0)})\otimes
I_{\cP_m}]\left[I_{\cE\otimes \cP_m}+\varphi(S_1^{(m)},\ldots,
S_n^{(m)})\right].
$$
where $\varphi(S_1^{(m)},\ldots, S_n^{(m)})=\sum_{k=1}^m
\sum_{|\alpha|=k} (I_\cE+B_{(0)})^{-1}B_{(\alpha)}\otimes
S_\alpha^{(m)}$. Since $[\varphi(S_1^{(m)},\ldots,
S_n^{(m)})]^{m+1}=0$, it is clear that the operator $I_{\cE\otimes
\cP_m}+\varphi(S_1^{(m)},\ldots, S_n^{(m)})$ is invertible and
$$
[I+f(S_1^{(m)},\ldots, S_n^{(m)})]^{-1}=[I_{\cE\otimes
\cP_m}+\varphi(S_1^{(m)},\ldots,
S_n^{(m)})]^{-1}[(I_\cE+B_{(0)})^{-1}\otimes I_{\cP_m}].
$$
Therefore, relation \eqref{Gm} implies
\begin{equation}
\label{a*} G_m^*=[I+f(S_1^{(m)},\ldots,
S_n^{(m)})^*]^{-1}[f(S_1^{(m)},\ldots, S_n^{(m)})^*-I].
\end{equation}

Now, notice that, for each $m\in \NN$ and $i=1,\ldots, n$, we have
$$(S_i^{(m+1)})^*|_{\cP^{(m)}}=(S_i^{(m)})^*=S_i^*|_{\cP^{(m)}}.
$$
Hence,
 $G_{m+1}^*|_{\cE\otimes\cP^{(m)}}=G_m^*$ for any $m\in \NN$.
 One can prove that there is a unique
contraction  $G\in  B(\cE\otimes F^2(H_n))$ such that
$G^*|_{\cE\otimes\cP^{(m)}}=G_m^*$ for any $m\in \NN$. Indeed, if
$x\in \cE\otimes F^2(H_n)$ let $q_m:=P_{\cE\otimes \cP^{(m)}} x$ and
notice that $\{G_m^* q_m\}_{m=1}^\infty$ is a Cauchy sequence.
Therefore, we can  define $G^*x:=\lim_{m\to \infty} G_m^* q_m$.
Since $\|G_m\|\leq 1$ for $m\in\NN$, so is the operator $G$.

Taking into account that $R_iS_j=S_jR_i$, $i,j=1,\ldots,n$, and that
$\cP^{(m)}$, $m\in\NN$, is an invariant subspace under each operator
$R_1^*,\ldots, R_n^*$, $S_1^*,\ldots, S_n^*$, we deduce that
$$(S_j^{(m+1)})^* R_i^*|_{\cP^{(m+1)}}=R_i^* (S_j^{(m+1)})^*\quad
\text{ for } \ i,j=1,\ldots,n.
$$
Hence and due to the form of  the operator $G_m$, we have
$$
(I_\cE\otimes R_i^*) G_m^*=G^*_m(I_\cE\otimes R_i^*) $$
 for any $m\in \NN$ and $i=1,\ldots, n$. Note that,
for each $\alpha\in\FF_n^+$ with $|\alpha|=k$, and $k=0,1,\ldots$,
we have
\begin{equation*}
\begin{split}
(I_\cE\otimes R_i^*)G^*(x\otimes e_{\alpha g_i})&= (I_\cE\otimes
R_i^*)G_{k+1}^* e_{\alpha g_i}=G_{k+1}^*
(I_\cE\otimes R_i^*)(x\otimes e_{\alpha g_i})\\
&=G_{k+1}^* (x\otimes e_\alpha) =G_k^*(x\otimes  e_\alpha)
\end{split}
\end{equation*}
and $ G^* (I_\cE\otimes R_i^*)(x\otimes e_{\alpha g_i})=
G^*(x\otimes e_\alpha) =G^*_k (x\otimes e_\alpha)$. Hence, we deduce
that
$$(I_\cE\otimes R_i^*)G^*(x\otimes e_{\alpha g_i})=G^*(I_\cE\otimes
R_i^*)(x\otimes  e_{\alpha g_i}), \quad i=1,\ldots,n.
$$
 On the
other hand, if $\alpha\in \FF_n^+$  has the form $g_{i_1}\cdots
g_{i_p}$ with $g_{i_p}\neq i$, then $G^* (I_\cE\otimes
R_i^*)(x\otimes e_{\alpha})=0$ and
$$(I_\cE\otimes R_i^*)G^*(x\otimes
e_{\alpha})=G_{k+1}^*(I_\cE\otimes  R_i^*)(x\otimes e_{\alpha })=0,
$$
which shows that $G^* (I_\cE\otimes R_i^*) (x\otimes
e_{\alpha})=(I_\cE\otimes R_i^*)G^*(x\otimes e_{\alpha})$.
Therefore,
$$
G(I_\cE\otimes R_i)=(I_\cE\otimes R_i) G,\quad i=1,\ldots,n.
$$
According to \cite {Po-analytic}, we deduce that $G$ is in
$B(\cE)\bar\otimes F_n^\infty$, the weakly closed algebra generated
by the spatial tensor product. Due to \cite{Po-holomorphic},
\cite{Po-pluriharmonic}, there is a unique $g\in
[H^\infty(B(\cH)^n_1)]_{\leq1}$ having the boundary function $G$,
i.e., $G=\text{\rm SOT-} \lim_{r\to 1} g(rS_1,\ldots, rS_n)$. Hence,
and using the fact that $G^*|_{\cE\otimes \cP^{(m)}}=G_m^*$, we
deduce that
\begin{equation*}
\begin{split}
G_m&=\text{\rm SOT-} \lim_{r\to 1}P_{\cP^{(m)}} g(rS_1,\ldots,
rS_n)|_{\cE\otimes \cP^{(m)}}\\
&=\lim_{r\to 1} g(rS_1^{(m)},\ldots, rS_n^{(m)})
=g(S_1^{(m)},\ldots, S_n^{(m)}).
\end{split}
\end{equation*}
Hence and using   relation  \eqref{a*}, we deduce that
\begin{equation*}
\begin{split}
g(S_1^{(m)},\ldots, S_n^{(m)})&=
 [f(S_1^{(m)},\ldots,
S_n^{(m)})-I][I+f(S_1^{(m)},\ldots, S_n^{(m)})]^{-1}\\
&=\Gamma(f)(S_1^{(m)},\ldots, S_n^{(m)})
\end{split}
\end{equation*}
for any $m\in \NN$. This  shows  that $\tilde
g=\widetilde\Gamma(\tilde f)$, where $\tilde f, \tilde g$  are the
power series associated with $f$ and $g$, respectively. On the other
hand, since $\tilde f\in \CC_{+}[Z_1,\ldots, Z_n]$, Proposition
\ref{Cayley-series} implies $\tilde g\in \CC_{-}[Z_1,\ldots, Z_n]$.
Consequently, we have $\Gamma(f)=g\in [H_-^\infty(B(\cH)^n_1)]_{\leq
1}$, which proves that the Cayley transform is well defined with
values in $[H_-^\infty(B(\cH)^n_1)]_{\leq 1}$.

To prove injectivity of $\Gamma$, let $f_1, f_2\in
Hol^+(B(\cH)^n_1)$ such that $\Gamma f_1= \Gamma f_2$.  Then
$\widetilde{\Gamma f_1}= \widetilde{\Gamma f_2}$ and due to
Proposition \ref{inverse}, we deduce that $\tilde f_1=\tilde f_2$,
which implies $f_1=f_2$.

Now we  prove that the noncommutative Cayley transform  is
surjective.  First, we consider the case when $g\in
H^\infty(B(\cH)^n_1)$ and $\|g\|<1$. Assume that  $g$   has  the
representation
$$g(X_1,\ldots, X_n):=\sum_{k=0}^\infty\sum_{|\alpha|=k} A_{(\alpha)}\otimes  X_\alpha,\quad
(X_1,\ldots, X_n)\in [B(\cH)^n]_1.
$$
Notice that $\|A_{(0)}\|<1$ and, consequently,  the operator
$I-A_{(0)}$ is invertible. Thus $ g\in
\left[H_-^\infty(B(\cH)^n_1)\right]_{\leq 1}$. Since $\|g\|<1$, it
is clear that $(1-g)^{-1} (1+g)\in B(\cE)\bar\otimes F_n^\infty$.
According to \cite{Po-analytic}, there is a unique Fourier
representation
\begin{equation}
\label{1-g}
 (1-g)^{-1} (1+g)=\sum_{\alpha\in\FF_n^+} B_{(\alpha)}
\otimes S_\alpha,
\end{equation}
where
$$ B_{(\alpha)}=P_{\cE\otimes \CC} (I_\cE\otimes S_\alpha^*)
(I-g )^{-1}(g +I)|_{\cE\otimes \CC},\quad \alpha\in \FF_n^+,
$$
with the property that the series $\sum_{k=0}^\infty
\sum_{|\alpha|=k} B_{(\alpha)}\otimes r^{|\alpha|} S_\alpha$ is
convergent in the operator norm topology  for each $r\in [0,1)$.
Therefore, the map
$$f(X_1,\ldots, X_n):=\sum_{k=0}^\infty
\sum_{|\alpha|=k} B_{(\alpha)}\otimes  X_\alpha, \qquad (X_1,\ldots,
X_n)\in [B(\cH)^n]_1, $$ is a free holomorphic function. Due to the
noncommutative von Neumann inequality \cite{Po-von}, we have
$\|g(rS_1,\ldots, rS_n)\|\leq \|g\|<1$ for any $r\in [0,1)$. Using
the functional calculus for row contractions \cite{Po-funct}, one
can easily see that
   the operator
$ [I-g(rS_1,\ldots, rS_n)]^{-1}[g(rS_1,\ldots, rS_n)+I]$ is in
$B(\cE)\otimes_{min}\cA_n$ and, due to \eqref{1-g}, it  has the
representation
 $\sum_{k=0}^\infty\sum_{|\alpha|=k} B_{(\alpha)}\otimes
r^{|\alpha|} S_\alpha $. Therefore,
\begin{equation}
\label{la}
 [I-g(rS_1,\ldots, rS_n)]^{-1}[g(rS_1,\ldots, rS_n)+I]
 =f(rS_1,\ldots,
rS_n)
\end{equation}
for any $ r\in[0,1)$,
  Now, we prove that
\begin{equation}
\label{g*g}
 f(X)^* +f(X)\geq 0 \quad \text{   for any  }
X\in [B(\cH)^n]_1.
\end{equation}
 Due to the properties of the noncommutative Poisson transform (see \cite{Po-poisson}),
  it is enough to show that
$$
f(rS_1,\ldots, rS_n)^* + f(rS_1,\ldots, rS_n)\geq 0\quad \text{ for
any  } \ r\in [0,1).
$$
Notice  first that relation \eqref{la} implies
  \begin{equation*}
  \begin{split}
  g(rS_1,\ldots, rS_n)[I+f(rS_1,\ldots, rS_n)]
  &=f(rS_1,\ldots, rS_n)-I.
\end{split}
  \end{equation*}
  Hence, we  have
  \begin{equation*}
  \begin{split}
2&[f(rS_1,\ldots, rS_n)^* + f(rS_1,\ldots,
rS_n)]\\
&=[I+f(rS_1,\ldots, rS_n)^*] [I+f(rS_1,\ldots,
rS_n)]-[I-f(rS_1,\ldots, rS_n)^*][I- f(rS_1,\ldots, rS_n)]\\
&=[I+f(rS_1,\ldots, rS_n)^*] [I+f(rS_1,\ldots, rS_n)]\\
 &\quad -[I+f(rS_1,\ldots, rS_n)^*] g(rS_1,\ldots, rS_n)^* g(rS_1,\ldots, rS_n)
  [I+f(rS_1,\ldots,
rS_n)]\\
&=[I+f(rS_1,\ldots, rS_n)^*][I-g(rS_1,\ldots, rS_n)^* g(rS_1,\ldots,
rS_n)] [I+f(rS_1,\ldots, rS_n)].
\end{split}
  \end{equation*}
Since $\|g(rS_1,\ldots, rS_n)\|< 1$, we deduce that $ f(rS_1,\ldots,
rS_n)^* + f(rS_1,\ldots, rS_n)\geq 0$  for any $r\in[0,1)$, which
proves inequality \eqref{g*g}.  This shows that $f\in
Hol^+(B(\cH)^n_1$ and, as in the proof of Proposition \ref{Re},  one
can see that the operator $I+f(rS_1,\ldots, rS_n)$ is invertible.
Now relation \eqref{la} implies
$$
g(rS_1,\ldots, rS_n)=[f(rS_1,\ldots, rS_n)-1][I+f(rS_1,\ldots,
rS_n)]^{-1}=(\Gamma f)(rS_1,\ldots, rS_n)
$$
for any $r\in [0,1)$. Hence and due to the uniqueness of Fourier
representation of the elements of $B(\cE)\bar \otimes F_n^\infty$,
we deduce that
    $\Gamma(f)=g$.

Now we consider the general case. Assume that $\psi\in
H_{-}^\infty(B(\cH)^n_1)$  has the representation
$\psi(X)=\sum_{k=0}^\infty \sum_{|\alpha|=k} B_{(\alpha)}\otimes
X_\alpha$ and $\|\psi\|\leq1$. For each $\epsilon\in [0,1)$ we  set
$g_\epsilon:=\epsilon \psi$.  Since $\|g_\epsilon\|<1$, we can apply
our previous result to $g_\epsilon$ and  find $f_\epsilon\in
Hol^+(B(\cH)^n_1)$ such that $\Gamma(f_\epsilon)=g_\epsilon $.
Assume that $f_\epsilon$ has the representation
$f_\epsilon(X):=\sum_{k=0}^\infty \sum_{|\alpha|=k}
A_{(\alpha)}^{(\epsilon)}\otimes X_\alpha$ for some
$A_{(\alpha)}^{(\epsilon)}\in B(\cE)$. Since  $\psi\in
[H_{-}^\infty(B(\cH)^n_1)]_{\leq 1}$ ,  the operator $I_\cE-B_{(0)}$
is invertible and $\|B_{(0)}\|\leq 1$. Taking into account that the
mapping \ $\epsilon\mapsto (I-\epsilon B_{(0)})^{-1}$\  is
continuous on $[0,1]$, we can find a constant $M>0$ such that
$\sup_{\epsilon\in [0,1]}\|(I-\epsilon B_{(0)})^{-1}\|\leq M$. On
the other hand, since  $g_\epsilon $ is the noncommutative Cayley
transform of $f_\epsilon$, we have $f_\epsilon(0)=(1-\epsilon
g(0))^{-1}(1+\epsilon g(0))$. Hence  we deduce that
$A_{(0)}^{(\epsilon)}=(1-\epsilon B_{(0)})^{-1}(1+\epsilon B_{(0)})$
and
$$
\left\|A_{(0)}^{(\epsilon)}\right\|\leq 2M \quad \text{ for any }\
\epsilon\in [0,1].
$$
On the other hand, since $\text{\rm Re}\, f_\epsilon\geq 0$, we can
use Proposition \ref{properties} and deduce that
\begin{equation}
\label{ine2} \left\|\sum_{|\alpha|=m}
\left(A_{(\alpha)}^{(\epsilon)} \right)^* A_{(\alpha)}
^{(\epsilon)}\right\|^{1/2}\leq
2\left\|A_{(0)}^{(\epsilon)}\right\|\leq 4M
\end{equation}
for any $m= 0, 1,\ldots $ and $\epsilon\in [0,1]$. Hence, we have
$\left\|A_{(\alpha)}^{(\epsilon)}\right\|\leq 4 M$ for any
$\epsilon\in[0,1]$ and $\alpha\in \FF_n^+$.

Due to Banach-Alaoglu theorem,  the ball $[B(\cE)]_{4M}^-$  is
compact
 in the $w^*$-topology. Since $\cE$ is a separable Hilbert space,
  $[B(\cE)]_{4M}^-$ is  a metric space in the $w^*$-topology which coincides with
   the weak operator topology on $[B(\cE)]_{4M}^-$. Consequently,
 the diagonal process guarantees the
existence  of a sequence $\{\epsilon_k\}_{k\geq 1} \subset [0,1)$
such that $\epsilon_k\to 1$ and, for each $\alpha\in \FF_n^+$,
$C_{(\alpha)}:=$WOT-$\lim\limits_{\epsilon_k\to
1}A_{(\alpha)}^{(\epsilon_k)} $ exists.
  The inequality \eqref{ine2} implies
 \begin{equation*}
 \left\|\sum_{|\alpha|=m} C_{(\alpha)}^*
C_{(\alpha)}\right\|^{1/2}\leq 4M
\end{equation*}
for any $m\geq 0$.  Consequently, we have $\limsup_{m\to\infty}
\left\|\sum_{|\alpha|=m} C_{(\alpha)}^*
C_{(\alpha)}\right\|^{1/2m}\leq 1$, which implies that
$$
g(X_1,\ldots, X_n):=\sum_{m=0}^\infty \sum_{|\alpha|=m}
C_{(\alpha)}\otimes X_\alpha,\quad (X_1,\ldots, X_n)\in
[B(\cH)^n]_1,
$$
is a free holomorphic function.

 Since  $\text{\rm Re}\, f_{\epsilon_k}(X_1,\ldots,
 X_n)\geq 0$ for any $(X_1,\ldots, X_n)\in [B(\cH)^n]_1$,
we have  $\text{\rm Re}\, f_{\epsilon_k}(S_1^{(p)},\ldots,
S_n^{(p)})\geq 0$ for any $p\geq 1$ and $k\geq 1$. Notice that, due
to the fact that $S_\alpha^{(p)}=0$ for any  $\alpha\in \FF_n^+$
with $|\alpha|\geq p+1$, we have $f_{\epsilon_k}(S_1^{(p)},\ldots,
S_n^{(p)})= \sum_{m=0}^p \sum_{|\alpha|=m}
A_{(\alpha)}^{(\epsilon_k)}\otimes S_\alpha^{(p)}$ and
 $g(S_1^{(p)},\ldots, S_n^{(p)})= \sum_{m=0}^p
\sum_{|\alpha|=m} C_{(\alpha)}\otimes S_\alpha^{(p)}$. Since for
each $\alpha\in \FF_n^+$,
$C_{(\alpha)}:=$WOT-$\lim\limits_{\epsilon_k\to
1}A_{(\alpha)}^{(\epsilon_k)}$, and $\text{\rm Re}\,
f_{\epsilon_k}(S_1^{(p)},\ldots, S_n^{(p)})\geq 0$,
 we deduce that
 $\text{\rm Re}\, g(S_1^{(p)},\ldots, S_n^{(p)})\geq 0$
for any $p\geq 1$.  Applying again Proposition \ref{properties}, we
deduce that
 $\text{\rm Re}\, g(X_1,\ldots,
 X_n)\geq 0$ for any $(X_1,\ldots, X_n)\in [B(\cH)^n]_1$.
Therefore, $g\in Hol^+(B(\cH)^n_1)$.

Now notice that
\begin{equation*}
\begin{split}
[\Gamma g](S_1^{(p)},\ldots, S_n^{(p)})&= \left[ g(S_1^{(p)},\ldots,
S_n^{(p)})-I\right]\left[I+g(S_1^{(p)},\ldots,
S_n^{(p)})\right]^{-1}\\
&= \text{\rm WOT-}\lim_{\epsilon_k\to 1}\left\{
 \left[ f_{\epsilon_k}(S_1^{(p)},\ldots,
S_n^{(p)})-I\right]\left[I+f_{\epsilon_k}(S_1^{(p)},\ldots,
S_n^{(p)})\right]^{-1}\right\}\\
&=\text{\rm WOT-}\lim_{\epsilon_k\to 1}[\Gamma
f_{\epsilon_k}](S_1^{(p)},\ldots, S_n^{(p)}) =\text{\rm
WOT-}\lim_{\epsilon_k\to 1} \epsilon_k
\psi(S_1^{(p)},\ldots, S_n^{(p)})\\
&=\psi(S_1^{(p)},\ldots, S_n^{(p)})
\end{split}
\end{equation*}
for any $p=0,1\ldots$. Hence, we deduce that $\Gamma g=\psi$. The
last part of the theorem is now obvious. The
 proof is complete.
\end{proof}

Let  $H_\CC^\infty(B(\cH)^n_1)$  be  the set of all bounded free
analytic functions with scalar coefficients, and denote by
$Hol_\CC^+(B(\cH)^n_1)$  the set of  all free holomorphic functions
$f$ with scalar coefficients such that $\text{Re}\, f\geq 0$.
Theorem \ref{Cayley1}, provides a noncommutative Cayley transform
$$
\Gamma_\CC:Hol_\CC^+(B(\cH)^n_1)\to \{f\in H_\CC^\infty(B(\cH)^n_1):
\|f\|\leq 1 \text{ and }\ f(0)\neq 1\}.
$$
In this case, we can say  more.
\begin{corollary} \label{scalar}
The Cayley transform $\Gamma_\CC$ can be identified with   a
bijection between $Hol_\CC^+(B(\cH)^n_1)$ and
$\left(F_n^\infty\right)_{\leq 1}\backslash\{I\}$, where
$F_n^\infty$ is the noncommutative analytic Toeplitz algebra.
\end{corollary}

\begin{proof}
According to \cite{Po-holomorphic}, the algebra
$H_\CC^\infty(B(\cH)^n_1)$ is completely isometric isomorphic  to
the noncommutative analytic Toeplitz algebra $F_n^\infty$.
Therefore,  due to Theorem \ref{Cayley1}, it is enough to show that
the set $\{f\in  F_n^\infty: \|f\|\leq 1 \text{ and }\ f(0)\neq 1\}$
coincides with $\left(F_n^\infty\right)_{\leq 1}\backslash\{I\}$.
Since one of the inclusions is obvious, assume that $g\in
\left(F_n^\infty\right)_{\leq 1}\backslash\{I\}$  and has the
representation $g=\sum_{\alpha\in \FF_n^+} a_\alpha S_\alpha$.
According to the Wiener type inequality for  the noncommutative
analytic Toeplitz algebra $F_n^\infty$  (see  \cite{PPoS}), we have
$$
\left(\sum_{|\alpha|=k} |a_\alpha|^2\right)^{1/2}\leq
1-|g(0)|^2\quad \text{ for any } \ k\geq 1.
$$
Consequently, if $g(0)=1$, then $a_\alpha=0$ for any $\alpha\in
\FF_n^+$ with $|\alpha|\geq 1$,  whence $g=I$, a contradiction.
Therefore, we must have $g(0)\neq 1$, which completes the proof.
\end{proof}

Denote by $\CC^{(m)}[Z_1,\ldots, Z_n]$, $m\in \NN$, the set
 of all noncommutative polynomials of degree $\leq m$.
 Let $\cA_{-}^{(m)}$ be the set of all operators
 $q(S_1^{(m)},\ldots, S_n^{(m)})\in B(\cE\otimes \cP^{(m)})$, where
  $q\in \CC^{(m)}[Z_1,\ldots, Z_n]$ and $I-q(0)$ is invertible. We also denote by $\cL^{(m)}$
the set of all operators $p(S_1^{(m)},\ldots, S_n^{(m)}) $ with
   $\text{\rm Re}\,p(S_1^{(m)},\ldots,
 S_n^{(m)})\geq 0$.
We consider  now the truncated (or constrained)   Cayley transforms
 $\cC^{(m)}$, $m\in \NN$,
 and point out the connection  with the noncommutative Cayley transform  $\Gamma$.

\begin{corollary}\label{Cayley2}
The Cayley transform $\cC^{(m)}:\cL^{(m)} \to [\cA_{-}^{(m)}]_{\leq
1}$ defined by
$$
\cC^{(m)} (Y):=(Y-I)(I+Y)^{-1}, \qquad Y\in \cL^{(m)},
$$
is a bijection and its inverse is given by
$$
[\cC^{(m)}]^{-1}(X):= (I+X)(I-X)^{-1}, \qquad X\in
[\cA_{-}^{(m)}]_{\leq 1}.
$$
Moreover, for each $m\in\NN$,
\begin{enumerate}
\item[(i)]\
$\cC^{(m)}[g(S_1^{(m)},\ldots, S_n^{(m)})]=[\Gamma g]
(S_1^{(m)},\ldots, S_n^{(m)})$ \ for any $g\in Hol^+(B(\cH)^n_1)$;
\item[(ii)]\
$ [\cC^{(m)}]^{-1}[f(S_1^{(m)},\ldots, S_n^{(m)})]=[\Gamma
f]^{-1}(S_1^{(m)},\ldots, S_n^{(m)}) $ \ for any $f\in
[H_-^\infty(B(\cH)^n_1)]_{\leq 1}$.
\end{enumerate}
\end{corollary}
\begin{proof}

If $Y\in \cL^{(m)}$, then $Y+Y^*\geq 0$ and, as in the proof of
Theorem \ref{Cayley1}, there exists a contraction $A_m:\cE\otimes
\cP^{(m)}\to \cE\otimes \cP^{(m)}$ such that $ A_m=(I-Y)(I+Y)^{-1}.
$
 It is easy to see that $A_m$ has the form
$p(S_1^{(m)},\ldots, S_n^{(m)})$ for some polynomial
  $p\in \CC^{(m)}[Z_1,\ldots, Z_n]$  with $I-p(0)$ invertible, which shows that
   $A_m\in[ \cA_{-}^{(m)}]_{\leq 1}$. Therefore $\cC^{(m)}$ is
   well-defined and $\cC^{(m)}(Y)=A_m$.
   As in the proof of Proposition \ref{Cayley-series}, one can prove that
   the Cayley transform  $\Gamma^{(m)}$  is one-to-one.

To prove the surjectivity of $\cC^{(m)}$, let $X=q(S_1^{(m)},\ldots,
S_n^{(m)})$ be in $ [\cA_{-}^{(m)}]_{\leq 1}$.
 Since  $I-q(0)$ invertible, so  is the operator     $I-X$   and
\begin{equation}\label{Y=}
Y:=(I+X)(I-X)^{-1}
\end{equation}
has the form $r(S_1^{(m)},\ldots, S_n^{(m)})$, where
  $r\in \CC^{(m)}[Z_1,\ldots, Z_n]$. Notice also that
  $X(I+Y)=Y-I$ and
  \begin{equation*}
  \begin{split}
  2(Y+Y^*)&=(I+Y)^*(I+Y)-(I-Y^*)(I-Y)\\
&=(I+Y)^*(I+Y)-(I+Y^*) X^*X(I+Y)\\
  &=(I+Y^*)(I-X^*X)(I+Y)\geq 0.
  \end{split}
  \end{equation*}
Therefore $Y\in \cL^{(m)}$ and relation  \eqref{Y=} implies
$X=(Y-I)(Y+I)^{-1}$. Therefore $\cC^{(m)}$ is a bijection.

If $g\in Hol^+(B(\cH)^n_1)$ then $g(S_1^{(m)},\ldots, S_n^{(m)})$ is
in $\cL_{n}^{(m)}$ and
\begin{equation*}
\begin{split}
(\Gamma g)(S_1^{(m)},\ldots, S_n^{(m)})&=[g(S_1^{(m)},\ldots,
S_n^{(m)})-I][I+g(S_1^{(m)},\ldots,
S_n^{(m)})]^{-1}\\
&=\cC^{(m)}[g(S_1^{(m)},\ldots, S_n^{(m)})].
\end{split}
\end{equation*}
Consequently,  item (i) holds.  Setting $f=\Gamma g$ in (i) and
using Theorem \ref{Cayley1}, one can easily  deduce item (ii).
    The proof is  complete.
\end{proof}

We mention that  the noncommutative Cayley transform   will be a key
tool in a forthcoming paper, where we study free pluriharmonic
majorants.

\bigskip

  \section{ Free holomorphic
  functions with positive real parts and Nevanlinna-Pick  interpolation }

In this section we  solve the  Nevanlinna-Pick  interpolation
problem for free holomorphic functions  with positive real parts on
the noncommutative  open  unit ball $[B(\cH)^n]_1$, and  obtain
several characterizations for the scalar representations of free
holomorphic functions with positive real parts.

We begin with a few preliminaries  concerning the eigenvectors of
the creation operators on the full Fock space.
 Let $\lambda:=(\lambda^{(1)},\ldots,
\lambda^{(n)})$ be in the open unit ball of $\CC^n$, i.e.,
 $$\BB_n:=\{z=(z_1,\ldots, z_n)\in \CC^n:\
\|z\|:=(|z_1|^2+\cdots + |z_n|^2)^{1/2}<1\}.$$  We define  the
vectors $z_\lambda\in F^2(H_n)$ by
$$
z_\lambda:=(1-\bar\lambda^{(1)} S_1-\cdots -
\bar\lambda^{(n)}S_n)(1)
$$
and  recall  (see \cite{Po-disc}) that
\begin{equation}
\label{eigen}
 S_i^* z_\lambda=\bar\lambda^{(i)} z_\lambda \quad
\text{ and } \quad R_i^* z_\lambda=\bar\lambda^{(i)} z_\lambda, \
\text{ for }\ i=1,\ldots,n,
\end{equation}
where $S_1,\ldots, S_n$ (resp. $R_1,\ldots, R_n$) are the left
(resp. right) creation operators on the full Fock space. If
$\varphi(X):=\sum_{k=0}^\infty \sum_{|\alpha|=k} A_{(\alpha)}\otimes
X_\alpha$ is a free holomorphic function on the operatorial ball
$[B(\cH)^n]_1$ with coefficients in $B(\cE)$, then
$$\varphi_r(S_1,\ldots, S_n):=\sum_{k=0}^\infty \sum_{|\alpha|=k}
A_{(\alpha)}\otimes r^{|\alpha|} S_\alpha$$
 is in $B(\cE)\otimes
\cA_n$ (the spatial tensor product),  for each $r\in [0,1)$.  One
can prove that
\begin{equation}\label{fir}
\left<\varphi_r(S_1,\ldots, S_n)^* (h\otimes z_\lambda ), h'\otimes
z_\mu \right>=\left< {\varphi_r(\lambda)}^*h, h'\right>\left<
z_\lambda, z_\mu\right>
\end{equation}
  for any $\lambda,\mu\in \BB_n$ and $h,h'\in \cK$, and
$$
\left<\varphi_r(\lambda)h, h'\right>=\left<\varphi_r(S_1,\ldots,
S_n) (h\otimes 1),h'\otimes z_\lambda
\right>=\left<\varphi_r(S_1,\ldots, S_n) (h \otimes u_\lambda
),h'\otimes u_\lambda \right>,
$$
where $u_\lambda:=\frac{z_\lambda}{\|z_\lambda\|}$.
 We remark that similar results  hold for
$\varphi_r(R_1,\ldots, R_n)$.

 In what follows we solve the   Nevanlinna-Pick
interpolation problem for free holomorphic functions  with positive
real parts on the noncommutative ball $[B(\cH)^n]_1$.
\pagebreak
\begin{theorem}
\label{Nev} Let $\lambda_0=0$, $\lambda_1,\ldots, \lambda_k$ be
distinct points in the open unit ball $\BB_n$  and let $W_0 $,
$W_1,\ldots, W_k$ be in $B(\cE)$, where $\cE$ is a separable Hilbert
space. Then there exists a free holomorphic function $f$ with
coefficients in $B(\cE)$ such that
\begin{enumerate}
\item[(i)]
$\text{\rm Re}\,f\geq 0$ and
\item[(ii)] $f(\lambda_j)=W_j$ for any $j=0,1,\ldots, k$,
\end{enumerate}
if and only if the operator matrix
\begin{equation}\label{WW}
\left[ \frac{W_i+ W_j^*}{1-\left< \lambda_i,
\lambda_j\right>}\right]_{(k+1)\times(k+1)}
\end{equation}
is positive semidefinite.
\end{theorem}

\begin{proof}

First, assume that $f$ is a free holomorphic function satisfying
properties (i) and (ii). Then, for any $r\in [0,1)$,
$f_r(S_1,\ldots, S_n)\in  B(\cE)\otimes \cA_n $  and, due to
relation \eqref{fir}, we have
$$
\left<f_r(S_1,\ldots, S_n)^*( h_i\otimes z_{\lambda_i} ), h_j\otimes
z_{\lambda_j}  \right>= \left<z_{\lambda_i},z_{\lambda_j}\right>
\left<f_r(\lambda_i)^*h_i, h_j\right>
$$
for any $h_i\in \cE$ and $i,j=0,1,\ldots, k$. Define  the vector
$\xi:=\sum_{j=0}^k h_j\otimes z_{\lambda_j} $, where  $h_j\in \cE$,
$j=0,1,\ldots, k$, and note that
\begin{equation*}
\begin{split}
\left<[f_r(S_1,\ldots, S_n)^*+f_r(S_1,\ldots, S_n)]\xi, \xi\right>&=
\sum_{i,j=1}^k \left<z_{\lambda_j},z_{\lambda_i}\right>
\left<(f_r(\lambda_i)+ f_r(\lambda_j)^*)h_j, h_i\right>\\
&=   \sum_{i,j=1}^k \left<\frac{ 1}{1-\left< \lambda_i,
\lambda_j\right>}(f_r(\lambda_i)+ f_r(\lambda_j)^*)h_j, h_i\right>.
\end{split}
\end{equation*}
On the other hand, since $\text{\rm Re}\, f\geq 0$, we have
$\left<[f_r(S_1,\ldots, S_n)^*+f_r(S_1,\ldots, S_n)]\xi,
\xi\right>\geq 0$ for any $r\in [0,1)$.  Consequently, the operator
matrix
$$
\left[ \frac{ f_r(\lambda_i)+ f_r(\lambda_j)^*}{1-\left< \lambda_i,
\lambda_j\right> }\right]_{(k+1)\times(k+1)}
$$
is positive semidefinite for each $r\in [0,1)$. Since $f$ is a free
holomorphic function, it is also  continuous   on $[B(\cH)^n]_1$.
Taking  $r\to 1$  we  deduce that
$$
\left[ \frac{ f(\lambda_i)+ f(\lambda_j)^*}{1-\left< \lambda_i,
\lambda_j\right> }\right]_{(k+1)\times(k+1)}
$$
is positive semidefinite.

Conversely, assume that the operator matrix  \eqref{WW} is positive
semidefinite. Consider  the  Hilbert space $\cN:=\text{\rm span}\,
\{z_{\lambda_j}: \ j=0,1,\ldots, k\}$. Since the vectors
$z_{\lambda_j}$,  $j=0,1,\ldots, k$, are linearly independent, one
can define the bounded linear operator $A:\cE\otimes \cN\to
\cE\otimes \cN$ by setting
\begin{equation}
\label{A*} A^*(h\otimes z_{\lambda_j} ):= W_j^*h\otimes
z_{\lambda_j}
\end{equation}
for any $h\in \cE$ and $j=0,1,\ldots, k$. Using again   the vector
$\xi:=\sum_{j=0}^k h_j\otimes z_{\lambda_j} $, where $h_j\in \cE$,
$j=0,1,\ldots, k$, we have
\begin{equation*}
\begin{split}
\left<(A^*+A)\xi, \xi\right>&= \left<\sum_{j=0}^k
 W_j^* h_j\otimes z_{\lambda_j}, \sum_{i=0}^k
  h_i\otimes z_{\lambda_i}\right>
  + \left<\sum_{j=0}^k h_j\otimes z_{\lambda_j}, \sum_{i=0}^k
 W_i^* h_i\otimes z_{\lambda_i}\right>\\
&= \sum_{i,j=0}^k \left<z_{\lambda_j},z_{\lambda_i}\right> \left<(
W_i+  W_j^*)h_j, h_i\right>\\
&= \sum_{i,j=0}^k   \left<\frac{ W_i+  W_j^*}{1-\left< \lambda_i,
\lambda_j\right>}h_j, h_i\right>.
\end{split}
\end{equation*}
 Consequently, $A^*+ A\geq 0$.  According to Proposition \ref{Re}, $I+A$  is invertible and the
 Cayley transform of $A$ is  the operator $T\in B(\cE\otimes \cN)$
 defined by
 $T=(A-I)(I+A)^{-1}$. Moreover, $T$ is a contraction with $I-T$
 invertible.
Since $W_j^*+W_j\geq 0$ for   $j=0,1,\ldots,n$, one can use again
Proposition \ref{Re}, to deduce that the operator $I+W_j^*$ is
invertible. Due to relation \eqref{A*}, we have
 $$
 (I+A^*)^{-1}(h\otimes z_{\lambda_j})=
 [(I+W_j^*)^{-1}h]\otimes z_{\lambda_j}
 $$
 and
 \begin{equation*}
 \begin{split}
T^*(h\otimes z_{\lambda_j})&= (A^*-I)(I+A^*)^{-1}
(h\otimes z_{\lambda_j})\\
&=(A^*-I)[ (I+W_j^*)^{-1}h\otimes z_{\lambda_j}]\\
&= [(W_j^*-I)(I+W_j^*)^{-1}h]\otimes z_{\lambda_j}
 \end{split}
 \end{equation*}
 for any $j=0,1,\ldots,k$.
 On the other hand,  the operator
$\Lambda_j^*:=(W_j^*-I)(I+W_j^*)^{-1}$ is the Cayley transform of
$W_j$ and,  consequently,  $\Lambda_j$ is a contraction with
$I-\Lambda_j$ invertible for each $j=0,1,\ldots, k$. Therefore, we
have
\begin{equation}
\label{T*}
 T^*(h\otimes z_{\lambda_j})=
\Lambda_j^*h\otimes z_{\lambda_j}, \quad j=0,1,\ldots, n.
\end{equation}
For each $i=1,\ldots, n$, define $X_i\in B(\cE\otimes \cN)$ by
$X_i:=I_\cE\otimes P_\cN S_i|_\cN $, where $S_1,\ldots, S_n$ are the
left creation operators on the full Fock space. Using relation
\eqref{eigen}, a simple calculation reveals that
$$
X_i^* T^*(h\otimes z_{\lambda_j})=T^*X_i^*(h\otimes z_{\lambda_j})
$$
for any $h\in \cE$ and $j=0,1,\ldots, k$. Since $\cN$ is an
invariant subspace under $S_1^*, \ldots, S_n^*$ and $T$ commutes
with each $X_1,\ldots, X_n$, one can apply the noncommutative
commutant lifting theorem \cite{Po-isometric} and use the
characterization of the commutant of $\{S_1,\ldots, S_n\}$ from
\cite{Po-analytic}, to find a contraction $\Phi(R_1,\ldots,
R_n):=\sum_{p=0}^\infty \sum_{|\alpha|=p}C_{(\alpha)}\otimes
R_\alpha$ in $B(\cE)\bar\otimes R_n^\infty$ such that
$$
P_{\cE\otimes \cN}\Phi(R_1,\ldots, R_n)|_{\cE\otimes \cN}=T,
$$
where $P_{\cE\otimes \cN}$ is the orthogonal projection onto
$\cE\otimes \cN$.  Since SOT-$\lim_{r\to 1}\Phi_r(R_1,\ldots,
R_n)=\Phi(R_1,\ldots, R_n)$ and using the remarks preceding the
theorem,
 we deduce that
\begin{equation*}
\begin{split}
\left<z_{\lambda_j}, z_{\lambda_j}\right> \left< \Phi(\lambda_j) h,
h'\right> &= \left<\Phi(R_1,\ldots, R_n)(h\otimes z_{\lambda_j} ),
(h'\otimes z_{\lambda_j})\right>\\
&=\left<T(h\otimes z_{\lambda_j}),
(h'\otimes z_{\lambda_j})\right>\\
&=\left<(h\otimes z_{\lambda_j}),
( \Lambda_j^*h'\otimes z_{\lambda_j})\right>\\
&=\left<z_{\lambda_j}, z_{\lambda_j}\right> \left< \Lambda_j h,
h'\right>
\end{split}
\end{equation*}
for any $h,h'\in \cK$ and $j=0,1,\ldots, k$. This clearly implies
$\Phi(\lambda_j)=\Lambda_j$ for any $j=0,1,\ldots, k$.

 Now, since
$I-\Phi(0)=I-\Lambda_0$ is invertible,  we can apply  Theorem
\ref{Cayley1} to $\Phi\in [H^\infty_-(B(\cH)^n_1]_{\leq 1}$,
  in order to find  a free holomorphic
function  $f$ with $\text{\rm Re}\, f\geq 0$  such that $\Gamma
f=\Phi$. Consequently, we have $\widetilde \Phi (1+\widetilde
f)=(\widetilde f-1)$, where $\widetilde\Phi$, $\widetilde f$ are the
power series associated  the free holomorphic functions $\Phi$ and
$f$, respectively. Since $\lambda_j\in \BB_n$, we deduce that
$$
\Phi(\lambda_j)[I+f(\lambda_j)]=f(\lambda_j)-I, \quad j=0,1,\ldots,
k.
$$
On the other hand,  taking into account that
$I-\Phi(\lambda_j)=I-\Lambda_j$ is invertible for each
$j=0,1,\ldots,n$, we deduce that
\begin{equation*}
\begin{split}
f(\lambda_j)&=[I-\Phi(\lambda_j)]^{-1}[\Phi(\lambda_j)+I)]\\
&=(I-\Lambda_j)^{-1} (\Lambda_j+I)=W_j.
\end{split}
\end{equation*}
This completes the proof.
\end{proof}

 An important connection to analytic function
theory in several complex variables is provided by the following
   characterizations  of  the scalar representations of free
holomorphic functions with positive real parts.
 \pagebreak
\begin{theorem}
\label{caracteriz} Let $\cE$ be   a separable Hilbert space and let
$F:\BB_n\to B(\cE)$  be  an operator-valued    function. Then the
following statements are equivalent:

\begin{enumerate}
\item[(i)]
 there exists a free holomorphic function $g$ with
 coefficients in $B(\cE)$
such that
$$
   \text{\rm Re}\, g \geq 0 \ \text{ and } \ F(z)=g(z),\
z\in \BB_n;
$$
\item[(ii)] there exists $\phi\in B(\cE) \bar \otimes F_n^\infty$ with
$\|\phi\|\leq 1$ and $I-\phi(0)$ invertible,  and such that
$$F(z)=[\Gamma^{-1}\phi](z)\quad \text{ for any } \  z\in \BB_n,
$$ where $\Gamma$ is the noncommutative Cayley transform.
\item[(iii)] $F$ is an operator-valued analytic  function on $\BB_n$ such
that
     the map
   \begin{equation}
   \label{FF*} \BB_n\times \BB_n\ni
(z,w)\mapsto \frac{F(z)+ {F(w)}^*}{1-\left<z,w\right>}\in B(\cE)
\end{equation}
 is positive  semidefinite;

\item[(iv)] there exists an $n$-tuple of isometries $(V_1,\ldots, V_n)$ on a
Hilbert space $\cK $
 with orthogonal ranges, and a bounded operator $W:\cE\to \cK$
such that
\begin{equation}
\label{F-rep}
 F(z_1,\ldots, z_n)=W^*(I+z_1V_1^*+\cdots
+z_nV_n^*)(I-z_1V_1^*-\cdots -z_nV_n^*)^{-1} W +i(\text{\rm Im}\,
F(0))
 \end{equation}
 for any $ (z_1,\ldots,
z_n)\in \BB_n$.
 \end{enumerate}
\end{theorem}

\begin{proof} The equivalence of  $(i)$ with  $(ii)$ is due to Theorem
\ref{Cayley1} and the identification of $H^\infty(B(\cH)^n_1)$ with
$B(\cE)\bar \otimes F_n^\infty$. The implication $(i)\implies (iii)$
follows as in the proof of Theorem \ref{Nev}.  Now, assume that $F$
satisfies condition $(iii)$ and let $\{\lambda_j\}_{j=0}^\infty$ be
a countable dense set in $\BB_n$ such that $\lambda_0:=0$. Since,
for each $k=0,1,\ldots$, the operator matrix
$$
\left[ \frac{ f(\lambda_i)+ f(\lambda_j)^*}{1-\left< \lambda_i,
\lambda_j\right> }\right]_{i,j=0,\ldots,k}
$$
is positive semidefinite, we can apply
  Theorem \ref{Nev}   to find a  free holomorphic function $f_k(X_1,\ldots,
X_n)=\sum_{m=0}^\infty \sum_{|\alpha|=m} A_{(\alpha)}^{(k)} \otimes
X_\alpha$ with coefficients $A_{(\alpha)}^{(k)}\in B(\cE)$ such that
\begin{equation}
\label{re=} \text{\rm Re}\, f_k\geq 0, \ \text{ and } \
f_k(\lambda_j)=F(\lambda_j)\quad \text{ for }\   j=0,1,\ldots, k.
\end{equation}
Let $C>0$ be such that $\|F(0)\|\leq C$. Applying  Proposition
\ref{properties} to $f_k$ , we deduce that
\begin{equation}
\label{ine22} \left\|\sum_{|\alpha|=m}
\left(A_{(\alpha)}^{(k)}\right)^*
A_{(\alpha)}^{(k)}\right\|^{1/2}\leq\|F(0)+F(0)^*\|\leq  2C
\end{equation}
for any $m\geq 0$. Hence, we have
$\left\|A_{(\alpha)}^{(k)}\right\|\leq 2C$ for any $k\geq 0$ and
$\alpha\in \FF_n^+$.

Due to Banach-Alaoglu theorem,  the ball $[B(\cE)]_{2C}^-$  is
compact
 in the $w^*$-topology. Since $\cE$ is a separable Hilbert space,
  $[B(\cE)]_{2C}^-$ is  a metric space in the $w^*$-topology which coincides with
   the weak operator topology on $[B(\cE)]_{2C}^-$. Consequently,
 the diagonal process guarantees the
existence  of a sequence $\{r_k\}\subset
 \NN$  such
that, for each  $\alpha\in \FF_n^+$ such that
$B_{(\alpha)}:=$WOT-$\lim\limits_{r_k\to \infty}A_{(\alpha)}^{(r_k)}
$ exists.
  The inequality \eqref{ine22} implies
 $
 \left\|\sum_{|\alpha|=m} B_{(\alpha)}^*
B_{(\alpha)}\right\|^{1/2}\leq 2C
$
for any $m\geq 0$.  Consequently,
$$
g(X_1,\ldots, X_n):=\sum_{m=0}^\infty \sum_{|\alpha|=m}
B_{(\alpha)}\otimes X_\alpha,\quad (X_1,\ldots, X_n)\in
[B(\cH)^n]_1,
$$
is a free holomorphic function.
 Since,  $\text{\rm Re}\, f_k(X_1,\ldots,
 X_n)\geq 0$ for any $(X_1,\ldots, X_n)\in [B(\cH)^n]_1$,
we have  $\text{\rm Re}\, f_k(S_1^{(p)},\ldots, S_n^{(p)})\geq 0$
for any $p\geq 1$ and $k\geq 0$.   Since, for each $\alpha\in
\FF_n^+$, $B_{(\alpha)}:=$WOT-$\lim\limits_{r_k\to
\infty}A_{(\alpha)}^{(r_k)}$,
 we deduce that
 $\text{\rm Re}\, g(S_1^{(p)},\ldots, S_n^{(p)})\geq 0$
for any $p\geq 1$. Due  to  Proposition \ref{properties}, we deduce
that
 $\text{\rm Re}\, g(X_1,\ldots,
 X_n)\geq 0$ for any $(X_1,\ldots, X_n)\in [B(\cH)^n]_1$.

 Let $\lambda=(\lambda^{(1)},\ldots \lambda^{(n)})\in \BB_n$ and
$r:=\|\lambda\|_2<1$. We use the notation
$\lambda_\alpha:=\lambda^{(i_1)}\cdots \lambda^{(i_k)}$  if
$\alpha=g_{i_1}\cdots g_{i_k}$ and $\lambda_{g_0}:=1$. Cauchy's
inequality and relation \eqref{ine22} imply
\begin{equation*}
\begin{split}
\left\|\sum_{|\alpha|=m} \lambda_\alpha A_{(\alpha)}^{(r_k)}
\right\|&\leq \left\|\sum_{|\alpha|=m}
\left(A_{(\alpha)}^{(r_k)}\right)^*
A_{(\alpha)}^{(r_k)}\right\|^{1/2}\left(\sum_{|\alpha|=m}|
\lambda_\alpha|^2\right)^{1/2}\\
&\leq 2C\left(\sum_{i=1}^n |\lambda^{(i)}|^2\right)^{m/2}=2Cr^m.
\end{split}
\end{equation*}
Similarly, we obtain $\left\|\sum_{|\alpha|=m} \lambda_\alpha
B_{(\alpha)} \right\|\leq 2Cr^m$. Consequently, for any $\epsilon>0$
there exists $N\in \NN$ such that
$$
\sum_{m\geq N+1}\left\|\sum_{|\alpha|=m} \lambda_\alpha
A_{(\alpha)}^{(r_k)} \right\|<\frac{\epsilon}{2} \quad \text{ and
}\quad \sum_{m\geq N+1}\left\|\sum_{|\alpha|=m}\lambda_\alpha
B_{(\alpha)} \right\|<\frac{\epsilon}{2}.
$$
Now, notice that
\begin{equation*}
\begin{split}
\|f_{r_k}(\lambda)-g(\lambda)\|&\leq \left\|\sum_{m=0}^N
\sum_{|\alpha|=m}\lambda_\alpha(A_{(\alpha)}^{(r_k)}-B_{(\alpha)})
\right\|+\epsilon\\
&\leq \sum_{m=0}^N
\left\|\sum_{|\alpha|=m}(A_{(\alpha)}^{(r_k)}-B_{(\alpha)})^*
(A_{(\alpha)}^{(r_k)}-B_{(\alpha)})\right\|^{1/2}\left(\sum_{|\alpha|=m}
|\lambda_\alpha|^2\right)^{1/2} +\epsilon \\
&\leq \sum_{m=0}^N
\left\|\sum_{|\alpha|=m}(A_{(\alpha)}^{(r_k)}-B_{(\alpha)})^*
(A_{(\alpha)}^{(r_k)}-B_{(\alpha)})\right\|^{1/2} +\epsilon.
\end{split}
\end{equation*}
 Since for each $\alpha\in \FF_n^+$,
$B_{(\alpha)}:=$WOT-$\lim\limits_{r_k\to
\infty}A_{(\alpha)}^{(r_k)}$,  we deduce that
$$
\text{\rm WOT-}\lim_{r_k\to \infty} f_{r_k}
(\lambda_j)=g(\lambda_j)\ \text{ for any }\ j=0,1,\ldots.
$$

Hence, and using \eqref{re=}, we obtain $g(\lambda_j)=F(\lambda_j)$
for any $j=0,1,\ldots$. Now, let us prove that $g(\zeta)=F(\zeta)$
for any $\zeta\in \BB_n$. To this end, let $\lambda\in \BB_n$ such
that $\lambda\neq \lambda_j$ for any $j=0,1,\ldots$. Applying the
preceding argument, we find a free holomorphic function $\varphi$
such that $\text{\rm Re}\, \varphi\geq 0$ and $\varphi(z)=F(z)$ for
any $z\in \{\lambda_j\}_{j=0}^\infty\cup \{\lambda\}$. Since, for
any $h,h'\in \cE$,  the maps $\zeta\mapsto
\left<\varphi(\zeta)h,h'\right>$ and $\zeta\mapsto \left<g(\zeta)h,
h'\right>$ are analytic on $\BB_n$, thus continuous,  and
$\{\lambda_j\}_{j=0}^\infty$ is dense in $\BB_n$, we must have
$\varphi(\zeta)=g(\zeta)$ for any $\zeta\in \BB_n$. In particular,
we have $g(\lambda)=\varphi(\lambda)=F(\lambda)$. Since $\lambda$ is
an arbitrary point in $\BB_n$, we deduce that $g(\lambda)
=F(\lambda)$ for any $\lambda\in \BB_n$. This completes the proof of
the implication $(iii)\implies (i)$.

Now, we prove that $(i)\implies (iv)$.  Since $g$ is a free
holomorphic function on $[B(\cH)^n]_1$  and $\text{\rm Re}\,g\geq
0$, Theorem 5.2 from \cite{Po-pluriharmonic}  shows that $g$ has
 the representation
 \begin{equation}
 \label{HR}
g(X_1,\ldots, X_n)=(\mu\otimes \text{\rm id})[ 2(I-R_1^*\otimes
X_1-\cdots -R_n^*\otimes X_n)^{-1}-I] +i(\text{\rm Im}\,
g(0))\otimes I
\end{equation}
for some  completely positive linear  map  $\mu$ from  the
Cuntz-Toeplitz algebra $C^*(R_1,\ldots, R_n)$ to $B(\cE)$. On the
other hand, due to Stinespring's representation theorem (see
\cite{St}),  there is a Hilbert space $\cK $, a bounded operator
$W:\cE\to \cK$, and
         a $*$-representation $\pi: C^*(R_1,\ldots, R_n)\to B(\cK)$    such that
         $$
         \mu(g)= W^*\pi(g)W  ,\quad  g\in C^*(R_1,\ldots, R_n).
        $$
        Notice that  $V_i:=\pi(R_i)$, $i=1,\ldots, n$, are isometries
        with orthogonal ranges.
Hence and using \eqref{HR}, we  obtain
 \begin{equation}
 \label{gProj}
g(X_1,\ldots, X_n)= (W^*\otimes I_\cH)[ 2(I-V_1^*\otimes X_1-\cdots
-V_n^*\otimes X_n)^{-1}-I]| (W\otimes I_\cH)+i(\text{\rm Im}\,
g(0))\otimes I.
\end{equation}
On the other hand,  since $V_1,\ldots, V_n$ are isometries with
orthogonal ranges and $(X_1,\ldots, X_n)\in [B(\cH)^n]_1$, we  have
$$ \|V_1^*\otimes X_1+\cdots +V_n^*\otimes
X_n\|=\left\|\sum_{i=1}^n X_iX_i^*\right\|^{1/2}<1.
$$
  Consequently,  we deduce that
 \begin{equation*}
 \begin{split}
  (I+V_1^*\otimes
X_1+\cdots +V_n^*\otimes X_n)&(I-V_1^*\otimes
X_1-\cdots -V_n^*\otimes X_n)^{-1}\\
&=[ 2(I-V_1^*\otimes X_1-\cdots -V_n^*\otimes X_n)^{-1}-I].
\end{split}
 \end{equation*}
Hence, and using   \eqref{gProj}, we obtain the representation
\eqref{F-rep}, if we take $\cH=\CC$. Therefore, $(iv)$ holds.

It remains to prove that $(iv)\implies (iii)$.  First, notice that
we can assume that $\text{\rm Im}\, F(0)=0$. Let $V_1,\ldots, V_n$
be  isometries with orthogonal ranges acting on a Hilbert space
$\cK$ and   $W:\cE\to \cK$  be a bounded operator such that,
\begin{equation*}
F(z_1,\ldots, z_n)=W^*(I+z_1V_1^*+\cdots
+z_nV_n^*)(I-z_1V_1^*-\cdots -z_nV_n^*)^{-1} W
 \end{equation*}
 for any $ (z_1,\ldots,
z_n)\in \BB_n$. Since  $\|z_1V_1^*+\cdots +z_nV_n^*\|<1$, we have
\begin{equation*}
\begin{split}
 (I_\cK+ z_1V_1^*+\cdots
+z_nV_n^*)&(I_\cK-z_1V_1^*-\cdots -z_nV_n^*)^{-1}\\
&=I_\cK+2\sum_{k=1}^\infty (z_1V_1^*+\cdots +z_nV_n^*)^k
=I_\cK+2\sum_{k=1}^\infty \left(\sum_{|\alpha|=k} z_\alpha
V_\alpha^*\right)
\end{split}
\end{equation*}
where the convergence  of the series is in the operator norm
topology of $B(\cK)$. Consequently,
$$
F(z_1,\ldots, z_n)= W^*W+2\sum_{k=1}^\infty \left(\sum_{|\alpha|=k}
z_\alpha  W^*V_\alpha^*W\right),\quad (z_1,\ldots, z_n)\in \BB_n.
$$
 This shows that $\BB_n\ni z\mapsto F(z_1,\ldots, z_n)\in B(\cE)$
is an operator-valued analytic function in $\BB_n$. On the other
hand, notice that
$$
G(X_1,\ldots, X_n):=2(I- S_1^*\otimes X_1-\cdots - S_n^*\otimes
X_n)^{-1}-I, \quad (X_1,\ldots, X_n)\in [B(\cH)^n]_1,
 $$
 is a free holomorphic function with coefficients in $B(F^2(H_n))$.
 Note also that
 \begin{equation*}
 \begin{split}
\text{\rm Re}\, G(X_1,\ldots, X_n)&=(I- S_1\otimes X_1^*-\cdots -
S_n\otimes X_n^*)^{-1}+(I- S_1^*\otimes X_1-\cdots -
S_n^*\otimes X_n)^{-1}-I\\
&=(I-S_1^*\otimes X_1-\cdots - S_n^*\otimes X_n)^{-1}[I\otimes
(I-X_1X_1^*-\cdots -X_nX_n^*]\\
&\qquad \qquad (I- S_1\otimes X_1^*-\cdots - S_n\otimes X_n^*)^{-1}
\end{split}
 \end{equation*}
and, therefore,  $\text{\rm Re}\, G(X_1,\ldots, X_n)\geq 0$ for any
$(X_1,\ldots, X_n)\in [B(\cH)^n]_1$. Now, we can apply the
implication $ (i)\implies (iii)$ and deduce that the map
\begin{equation}
   \label{GG*} \BB_n\times \BB_n\ni
(z,w)\mapsto \frac{G(z)+ {G(w)}^*}{1-\left<z,w\right>}\in B(\cE)
\end{equation}
 is positive  semidefinite.  On the other hand, since $V_1,\ldots,
 V_n$ are isometries with orthogonal ranges, there exists a unique
 $*$-representation $\pi:C^*(S_1,\ldots, S_n)\to B(\cK)$ such that
 $\pi(S_\alpha S_\beta^*)=V_\alpha V_\beta^*$ for any
 $\alpha,\beta\in \FF_n^+$ (see \cite{Cu}). Consider the completely positive linear
 map $\omega_W:B(\cK)\to B(\cE)$ defined  by
 $\omega_W(Y):=W^* Y W$. Using the above calculations, we
 deduce that
 \begin{equation*}
 \begin{split}
F(z) &= \omega_W\left[ I_\cK+2\sum_{k=1}^\infty
\left(\sum_{|\alpha|=k}
z_\alpha V_\alpha^*\right)\right]\\
&= (\omega_W\circ \pi)\left[ I_{F^2(H_n)}+2\sum_{k=1}^\infty
\left(\sum_{|\alpha|=k} z_\alpha S_\alpha^*\right)\right]\\
&=(\omega_W\circ \pi)(G(z))
\end{split}
 \end{equation*}
for any $z=(z_1,\ldots, z_n)\in \BB_n$. Hence, and using the fact
that the map \eqref{GG*} is positive semidefinite and $\omega_W\circ
\pi$ is a completely positive linear map, we deduce that the map
\eqref{FF*} is positive semidefinite. The proof is complete.
\end{proof}

We remark that if $F(0)=I$ in Theorem \ref{caracteriz}, then the
representation  \eqref{F-rep} becomes
$$ F(z_1,\ldots, z_n)=P_\cE(I+z_1V_1^*+\cdots
+z_nV_n^*)(I-z_1V_1^*-\cdots -z_nV_n^*)^{-1}|_\cE
 $$
 for any $ (z_1,\ldots,
z_n)\in \BB_n$, where  $\cK\supset \cE$ and $P_\cE$ is the
orthogonal projection of $\cK$ onto $\cE$.

A closer look at the proof of Theorem \ref{caracteriz} and  using
Corollary \ref{scalar}, enable us to deduce the following scalar
version of Theorem \ref{caracteriz}.

\begin{corollary}
\label{scalar2} Let $f:\BB_n\to \CC$ be a complex-valued  function.
Then the following statements are equivalent:

\begin{enumerate}
\item[(i)]
 there exists a free holomorphic function $g$ with  scalar
 coefficients
such that
$$
  \text{\rm Re}\, g \geq 0 \ \text{ and } \ g(z)=f(z),\
z\in \BB_n;
$$
\item[(ii)] there exists $\phi\in F_n^\infty$ with
$\|\phi\|\leq 1$, $\phi\neq I$, and such that
$f(z)=[\Gamma^{-1}\phi](z),\ z\in \BB_n$, where $\Gamma$ is the
noncommutative Cayley transform;
\item[(iii)]
  $f$ is analytic  and
 has the property that   the map $$ \BB_n\times \BB_n\ni
(z,w)\mapsto \frac{f(z)+ \overline{f(w)}}{1-\left<z,w\right>}\in \CC
$$
 is positive  semidefinite;

\item[(iv)] there exists an $n$-tuple of isometries $(V_1,\ldots, V_n)$ with orthogonal ranges on a
Hilbert space $\cK$,  a vector $\xi\in \cK$ such that
\begin{equation*}
 f(z )=\left<(I+z_1V_1^*+\cdots +
z_nV_n^*)(I-z_1V_1^*-\cdots - z_nV_n^*)^{-1}\xi, \xi\right> +
i(\text{\rm Im}\, f(0))
 \end{equation*}
 for  any $z:=(z_1,\ldots, z_n)\in \BB_n$;

\item[(v)] there exists  a completely positive linear map
$\nu:C^*(R_1,\ldots,R_n)\to \CC$   such that
$$
f(z)=\nu[(I+z_1R_1^*+\cdots + z_nR_n^*)(I-z_1R_1^*-\cdots -
z_nR_n^*)^{-1}]+  i(\text{\rm Im}\, f(0))
$$
for  any $z:=(z_1,\ldots, z_n)\in \BB_n$.
 \end{enumerate}
\end{corollary}

We mention that the implication $(iii)\implies (iv)$ was also proved
in \cite{McP} using different techniques.

\bigskip

\section{Free holomorphic (resp.~analytic) extensions}

We need to recall (see  \cite{Ru}) a few facts concerning the
automorphisms of the unit ball  $\BB_n$. Let $a\in \BB_n$ and
consider $\psi_a\in
 Aut(\BB_n)$, the automorphism of the unit ball, defined by
 \begin{equation}
 \label{auto}
 \psi_a(z):= \frac{a-P_az-s_a (I-P_a)z}{1-\left<z,a\right>}, \quad z\in
 \BB_n,
 \end{equation}
where   $P_0=0$, $P_a z:=\frac{\left<z,a\right>}{\left<a,a\right>}
a$ if $ a\neq 0$, and $s_a:=(1- \left<a,a\right>)^{1/2}$. The
automorphism $\psi_a:\BB_n\to \BB_n$ has the following properties:
\begin{enumerate}
\item[(i)]
$\psi_a(0)=a$, $\psi_a(a)=0$;
\item[(ii)] $\psi_a(\psi_a(z))=z$ for any $z\in \BB_n$;
\item[(iii)] $1-\left<\psi_a(z), \psi_a(w)\right>=
\frac{(1- \left<a,a\right>)(1- \left<z,w\right>)}{(1-
\left<z,a\right>)(1- \left<a,w\right>)} $ for any $z,w\in \BB_n$.
\end{enumerate}

  The main result    of this paper   is  the following  theorem
  regarding
   free holomorphic (resp.~analytic) extensions of
operator-valued functions defined on subsets of $\BB_n$.

\pagebreak

\begin{theorem}
\label{Lam} Let $\Lambda$ be an arbitrary  subset of $\BB_n$ and let
$F:\Lambda\to B(\cE)$ be an operator-valued function, where $\cE$ is
a separable Hilbert space. Then the following statements are
equivalent:

\begin{enumerate}
\item[(i)] $F$ has
   a free holomorphic  extension  $G$  with positive real part,
   i.e.,
$$
   \text{\rm Re}\, G \geq 0 \ \text{ and } \ G(z)=F(z),\
z\in \Lambda;
$$
\item[(ii)]
    the map
   \begin{equation*}
    \Lambda\times \Lambda \ni
(z,w)\mapsto \frac{F(z)+ {F(w)}^*}{1-\left<z,w\right>}\in B(\cE)
\end{equation*}
 is positive  semidefinite;

\item[(iii)] $F$ has an operator-valued  analytic extension
     $\Phi:\BB_n\to B(\cE)$ such that
    the map
  $$ \BB_n\times \BB_n\ni
(z,w)\mapsto \frac{\Phi(z)+ {\Phi(w)}^*}{1-\left<z,w\right>}\in
B(\cE)
$$
 is positive  semidefinite.
\end{enumerate}
\end{theorem}
\begin{proof}  Consider the case when  $0\in \Lambda$.
 If $\Lambda\subset \BB_n$ is finite, then the
result follows from Theorem \ref{Nev} and Theorem \ref{caracteriz}.
Now, assume that  $\Lambda$ is infinite. If condition  (ii) holds,
we choose a sequence $\{\lambda_j\}_{j=0}^\infty\subset \Lambda$
dense in $\Lambda$ and such that $\lambda_0=0$. As in the proof of
Theorem \ref{caracteriz}, one can find a free holomorphic function
$F$ with $\text{\rm Re}\, F\geq 0$ and $G(z)=F(z)$ for any $z\in
\Lambda$.  Therefore, item (i) holds.  The    other implications can
be proved similarly to those of Theorem \ref{caracteriz}.

 Consider the case when $\Lambda$  is an arbitrary  subset of
 $\BB_n$.  We prove the implication $(ii)\implies (i)$. Fix a point $a\in \Lambda$ and
 let  $\psi_a\in
 Aut(\BB_n)$  be the automorphism of the unit ball defined by relation \eqref{auto}.
Consider $\Lambda_0:=\psi_a(\Lambda)\subset \BB_n$ and define the
function $\Phi:\Lambda_0\to B(\cE)$ by setting
\begin{equation}
\label{Phi-xi}
 \Phi(\xi):=F(\psi^{-1}_a(\xi)), \quad \xi\in \Lambda_0.
\end{equation}
Assume that condition (ii) holds.  First, we prove that the map
\begin{equation}
\label{Lam0}
 \Lambda_0\times \Lambda_0 \ni (\xi,\eta)\mapsto
\frac{\Phi(\xi)+
 {\Phi(\eta)}^*}{1-\left<\xi,\eta\right>}\in B(\cE)
\end{equation}
 is positive  semidefinite. To this end, let $\xi_1,\ldots, \xi_k\in
 \Lambda_0$ and $h_1,\ldots, h_k\in \cE$, and
 set $z_i:=\psi_a^{-1}(\xi_i)\in \Lambda$ for $i=1,\ldots,k$.
  Using the properties of the automorphism $\psi_a$, preceding this theorem, we deduce that
 \begin{equation*}
 \begin{split}
\sum_{i=1}^k \sum_{j=1}^k \left<\frac{\Phi(\xi_i)+
 {\Phi(\xi_j)}^*}{1-\left<\xi_i,\xi_j\right>} h_i,h_j\right>
 &=\sum_{i=1}^k \sum_{j=1}^k \left<
 \frac{F(z_i)+
 {F(z_j)}^*}{1-\left<z_i,z_j\right>}\frac{1-\left<z_i,z_j\right>}{1-
 \left<\psi_a(z_i),\psi_a(z_j)\right>}h_i,h_j\right>\\
 &=
 \sum_{i=1}^k \sum_{j=1}^k \left<
 \frac{F(z_i)+
 {F(z_j)}^*}{1-\left<z_i,z_j\right>}\frac{(1-\left<z_i,a\right>)(1-\left<a,z_j\right>)}
 {1-\left<a,a\right>} h_i,h_j\right>\\
&=\sum_{i=1}^k \sum_{j=1}^k \left<
 \frac{F(z_i)+
 {F(z_j)}^*}{1-\left<z_i,z_j\right>}x_i,x_j\right>\geq 0,
\end{split}
 \end{equation*}
where
$$
x_i:=\frac{1-\left<z_i,a\right>}{\sqrt{1-\left<a,a\right>}} h_i\quad
\text{ for } \ i=1,\ldots,k.
$$
This proves our assertion.

Now, since $0=\psi_a(a)\in \Lambda_0$ and the map given by
\eqref{Lam0} is positive semidefinite, we can apply the first part
of the proof to find an operator-valued analytic  function
$\widetilde \Phi:\BB_n\to B(\cE)$   such that the map

\begin{equation}
\label{tildaPhi}
\BB_n\times \BB_n \ni (\xi,\eta)\mapsto
\frac{\widetilde\Phi(\xi)+
 {\widetilde\Phi(\eta)}^*}{1-\left<\xi,\eta\right>}\in B(\cE)
\end{equation}
 is positive  semidefinite and such that $\widetilde \Phi$ extends
 $\Phi$, i.e.,
 \begin{equation}
 \label{ext}
\widetilde \Phi(\xi)= \Phi(\xi)\quad \text{ for any }\ \xi\in
\Lambda_0.
 \end{equation}
Define the operator-valued analytic function  $g:\BB_n\to B(\cE)$ by
setting
\begin{equation}
\label{gz} g(z):=\widetilde \Phi(\psi_a(z)),\quad z\in \BB_n.
\end{equation}
Let $z_1,\ldots, z_k\in
 \BB_n$ and $h_1,\ldots, h_k\in \cE$, and
 set $\eta_i:=\psi_a^{-1}(z_i)$ for $i=1,\ldots,k$.
  Using the properties of the automorphism $\psi_a$, similar
  calculations as above imply
$$
\sum_{i=1}^k \sum_{j=1}^k \left<\frac{g(z_i)+
 {g(z_j)}^*}{1-\left<z_i,z_j\right>} h_i,h_j\right>
 =
 \sum_{i=1}^k \sum_{j=1}^k \left<\frac{\widetilde\Phi(\eta_i)+
 {\widetilde\Phi(\eta_j)}^*}{1-\left<\eta_i,\eta_j\right>} y_i,y_j\right>
 $$
 where
$$
y_i:=\frac{1-\left<\eta_i,a\right>}{\sqrt{1-\left<a,a\right>}}
h_i\quad \text{ for } \ i=1,\ldots,k.
$$
Since the map given by \eqref{tildaPhi} is positive semidefinite, so
is the map
$$
\BB_n\times\BB_n\ni(z,w)\mapsto
\frac{g(z)+g(w)^*}{1-\left<z,w\right>}\in B(\cE).
$$
  Applying  Theorem \ref{caracteriz} to the analytic function $g$,
we find a free holomorphic function $G$ with coefficients in
$B(\cE)$ with $\text{\rm Re}\, G\geq 0$ and $G(z)=g(z)$ for any
$z\in \BB_n$. Hence and using relations \eqref{Phi-xi}, \eqref{ext},
and  \eqref{gz}, we deduce that
$$
G(z)=g(z)=\widetilde\Phi(\psi_a(z))=\Phi(\psi_a(z))=F(z)
$$
for any $z\in \Lambda$, which proves the implication $(ii)\implies
(i)$. The implication  $(i)\implies (ii)$ follows as in the proof of
Theorem \ref{Nev}. Since  $(iii)\implies (ii)$ is obvious and
implication  $(i)\implies (iii)$ was proved in Theorem
\ref{caracteriz}, the proof is complete.
\end{proof}

\bigskip

       %

      \end{document}